\documentclass[12pt,a4paper,reqno]{amsart}
\usepackage[T1]{fontenc}
\usepackage[utf8]{inputenc}
\IfFileExists{newtxtext.sty}{\usepackage{newtxtext,newtxmath}}{\usepackage{lmodern}}
\usepackage{amsmath,amsthm,mathtools}
\usepackage[a4paper,left=2.4cm,right=2.4cm,top=2.2cm,bottom=2.0cm]{geometry}
\usepackage[colorlinks=true,allcolors=black]{hyperref}
\usepackage{microtype}
\usepackage{needspace}
\newtheorem{theorem}{Theorem}[section]
\newtheorem{proposition}[theorem]{Proposition}
\newtheorem{lemma}[theorem]{Lemma}
\newtheorem{corollary}[theorem]{Corollary}
\theoremstyle{remark}
\newtheorem{remark}[theorem]{Remark}
\newcommand{\OO}{\mathcal O}
\newcommand{\FF}{\mathcal F}
\newcommand{\CC}{\mathbb C}
\newcommand{\mm}{\mathfrak m}
\newcommand{\ord}{\operatorname{ord}}
\newcommand{\GSV}{\operatorname{GSV}}
\title[FFTS]{Foliations and Free Tangent Sheaves}
\author{Maur\'icio Corr\^ea}
\thanks{Dipartimento di Matematica, Universit\`a degli Studi di Bari Aldo Moro, Via Edoardo Orabona 4, 70125 Bari, Italy. E-mail: \texttt{mauricio.barros@uniba.it}.}
\date{30 November 2026}
\newtheorem*{maintheorema}{Theorem A}
\newtheorem*{maintheoremb}{Theorem B}
\begin{document}
\begin{abstract}
We study the Lipman--Zariski problem for normal complex surface germs with
locally free tangent sheaf through the foliations induced by their derivations.
We prove that, at a singular point, the saturated foliation induced by any
derivation which can be completed to a local frame of the tangent sheaf has a
dicritical component on the minimal good resolution. In particular, such a
germ is smooth if it admits such a derivation whose foliation has no dicritical
exceptional component or only finitely many invariant irreducible curve germs
through the point. We relate this obstruction to exceptional cycles and local
indices, determine the generic dicritical components for vector fields with
coefficients of prescribed order, and establish a Baum--Bott--Camacho--Sad
identity on compact normal surfaces.
\end{abstract}
\maketitle
\markboth{Foliations and Free Tangent Sheaves}{Foliations and Free Tangent Sheaves}
\section{Introduction}\label{sec:intro}
The Lipman--Zariski problem asks whether a normal complex space with
locally free tangent sheaf must be smooth. For a normal surface germ this
question has a natural formulation in terms of foliations. Local tangent
fields lift to logarithmic vector fields on the minimal good resolution
\cite[Fact~1]{BergnerGraf20}. Their logarithmic lifts may vanish along the
exceptional curves, and removing the divisorial factors as ordinary vector
fields can change the behaviour of the induced foliation: an exceptional
component may be invariant or generically transverse to it.

The conjecture has been proved under several hypotheses. Greb, Kebekus,
Kov\'acs and Peternell established it for klt spaces
\cite{GrebKebekusKovacsPeternell11}; Druel and, independently, Graf and
Kov\'acs treated log canonical spaces \cite{Druel14,GrafKovacs14}.
K\"allstr\"om proved the case of local complete intersections
\cite{Kallstrom11}. For algebraic surfaces, Biswas, Gurjar and Kolte
proved the conjecture for affine surfaces whose smooth locus has
logarithmic Kodaira dimension at most one, and for several classes of
projective surfaces \cite{BiswasGurjarKolte14}. Ballico excluded singular normal projective surfaces with locally free,
weakly ample tangent sheaf \cite{Ballico06}. Graf settled the cases of surface singularities
of genus one and of genus two whose links are not rational homology spheres
\cite{Graf21}. Bergner and Graf extended this to $p_g-g-b\leq 2$,
where $p_g$ is the geometric genus, $g$ the sum of exceptional genera,
and $b$ the first Betti number of the dual graph \cite{BergnerGraf20}.
Graf and Kov\'acs also proved the conjecture for potentially Du Bois
singularities \cite{GrafKovacsDuBois14}. For logarithmic
tangent sheaves on pairs, Bergner established toroidal and smoothness
results under log canonical and dlt hypotheses, respectively
\cite{Bergner21}.

In joint work with Louren\c{c}o and Machado
\cite{CorreaLourencoMachado26}, we studied logarithmic Baum--Bott residues
on singular varieties, their relation to the GSV and Camacho--Sad indices,
and a foliated smoothness criterion for compact surfaces. The subsequent
study of excess logarithmic residues with Cavalcante, Louren\c{c}o and
Shahsavaripour also relates exceptional residues to log discrepancies
\cite[Theorem~8.1]{CavalcanteCorreaLourencoShahsavaripour26}.
Here we study the local geometry of the lifted tangent fields themselves,
related to J\"order's criterion involving commuting tangent frames
\cite{Joerder14} and the geometry of separatrices on singular surfaces
\cite{Santos20}.

A preprint of OpenAI, dated 23 September 2026, proposes a singular normal
affine complex surface with free tangent sheaf
\cite{OAILipmanZariski}. If the proposed construction is valid, it would
provide a counterexample to the unrestricted Lipman--Zariski conjecture in
characteristic zero. It is therefore natural to ask which properties of the
foliations induced by a free tangent frame remain constrained by the
singular geometry. Our general results apply to arbitrary normal surface
germs with locally free tangent sheaf. The later realisation of prescribed
indices on the proposed surface is stated under the specific analytic and
algebraic hypotheses of that construction.

Let $(Y,p)$ be a normal complex surface germ with locally free tangent
sheaf $T_Y$, and let $\pi:X\to Y$ be its minimal good resolution, with
reduced exceptional divisor $E$. A vector field $v$ is \emph{primitive} if
it can be completed to a local frame of $T_Y$. Equivalently, when
$v=fe_1+ge_2$ in a local frame, at least one of $f,g$ is a unit at $p$,
or $v\notin\mathfrak m_pT_{Y,p}$. An exceptional component is
\emph{dicritical} if the ordinary saturated foliation of the lifted field
is generically transverse to it. Our principal result is the following.

\begin{maintheorema}
If $(Y,p)$ is singular and $T_Y$ is locally free, the saturated foliation of every primitive derivation has a dicritical component on the minimal good resolution. For every local frame $(e_1,e_2)$, there is a component $E_j$ which is dicritical for $e_1+\lambda e_2$ for all but finitely many $\lambda\in\mathbb C$.
\end{maintheorema}

In particular, the Lipman--Zariski conjecture holds for normal surface
germs with locally free tangent sheaf admitting a primitive derivation
whose saturated foliation has no dicritical exceptional component.
The existence of a primitive derivation with only finitely many invariant
irreducible analytic curve germs through $p$ likewise forces the surface
to be smooth (Corollary~\ref{cor:separatrices}).

The proof compares the logarithmic zero divisor of a lifted field with
the determinant divisor of a lifted tangent frame. If every exceptional
component is invariant, restriction to the exceptional curves gives an
effective anti-nef correction cycle $B_v=Z_v-W$. Its vanishing would
produce a logarithmic one-form with nonzero exceptional residue,
contradicting the negative definiteness of the intersection matrix.
The classical \emph{fundamental cycle} of Artin, $F_\pi=\sum_i\varphi_iE_i$,
is the coefficientwise smallest nonzero effective exceptional cycle with
$F_\pi\cdot E_i\leq0$ for every exceptional component
\cite{Artin66,Laufer71}; see also \cite[\S2.4]{NagyNemethi19}.
Consequently $B_v\geq F_\pi$, giving explicit lower bounds for the
coefficients of the divisorial zero cycle and of the coefficient ideal
(Corollary~\ref{cor:fundamental-bound}). A primitive field at a singular
point consequently has infinitely many invariant irreducible curve germs
through that point, obtained from leaves transverse to a dicritical
exceptional component.

For nonprimitive fields we compare the valuations of the coefficient
ideal with the Smith exponents of the logarithmic lifting matrix and the
maximal-ideal cycle. For general coefficients of prescribed order $k$,
this determines the logarithmic zero divisor exactly; the exceptional
components that are dicritical for a general field are independent of
$k$. The GSV indices on invariant curves and the classical tangency numbers
on dicritical curves are given by intersection numbers on the resolution. The distinction between the smooth-branch GSV index and the
Camacho--Sad index at an invariant crossing is essential in these
formulae \cite{GSVref,CamachoSad82,Brunella97}.

\Needspace{10\baselineskip}
Primitivity is essential to Theorem~A. A first integral on a normal
Gorenstein surface singularity leads to a derivation whose saturated
foliation preserves the exceptional divisor. The construction gives a
more precise complementary result.

\begin{maintheoremb}
Let $(Y,p)$ be an isolated normal Gorenstein algebraic surface singularity and $\pi:X\to Y$ a prescribed projective log resolution. There exists a regular derivation $v$ on a neighbourhood of $p$, with no zero on the punctured germ, whose saturated foliation has a holomorphic first integral, preserves every exceptional component, and has only nondegenerate singularities with negative rational Camacho--Sad indices above $p$. If $T_Y$ is locally free, this derivation is nonprimitive.
\end{maintheoremb}

We also investigate the realisation of prescribed exceptional
Camacho--Sad indices using nonprimitive fields in the analytic model of
\cite[Proposition~5.4]{OAILipmanZariski}, their behaviour under further
blow-ups, and their algebraic realisation subject to the additional
hypotheses of \cite[Proposition~6.2]{OAILipmanZariski}. This distinguishes
the general first-integral construction from the construction using
specific first-jet data of the proposed model.

On compact normal surfaces the same exceptional correction enters a
Baum--Bott--Camacho--Sad identity relating characteristic numbers and
exceptional indices. Unlike the general logarithmic residue comparison
of \cite[Corollaries~1.6 and 1.7]{CorreaLourencoMachado26}, the formula
here records the contribution of the cycle $B_v$ explicitly. For
negative real exceptional indices it becomes a sum of squares whose
equality case forces prescribed rational indices. The resulting
smoothness criterion is a numerical consequence under additional
hypotheses, and should be read alongside the generalised-curve criterion
of \cite[Corollary~1.9]{CorreaLourencoMachado26}. The comparison of
ordinary and logarithmic normal sheaves in
\cite[Theorem~1.1]{CavalcanteCorreaLourencoShahsavaripour26} gives further
context for this identity; the point here is the explicit correction
cycle and its quadratic equality criterion.

\section{Vector fields, foliations and local indices}\label{sec:prelim}

Let $Y$ be a normal complex surface. Its tangent sheaf is the sheaf of
holomorphic derivations $T_Y=\operatorname{Der}_{\mathbb C}(\OO_Y)$; a vector
field on $Y$ is a section of this sheaf. On the smooth locus, a nonzero field
$v$ determines a one-dimensional foliation after its common divisorial factors
have been removed. More generally, on a smooth surface $X$ a foliation by
curves is described by a saturated rank-one subsheaf
$T_{\FF}\subset T_X$. It is a line bundle, and in local coordinates $(x,y)$
we may choose a generator
$\xi=A(x,y)\partial_x+B(x,y)\partial_y$, where $A$ and $B$ have no
common irreducible factor. The associated one-form is
$\omega=B\,dx-A\,dy$; its singularities are the isolated common zeroes
of $A$ and $B$. Multiplication of $\xi$ by a unit does not change the
foliation.

Let $\pi\colon X\to Y$ be a good resolution of a normal surface germ,
with reduced simple-normal-crossings exceptional divisor $E$. A field on $Y$
may admit a holomorphic lift tangent to $E$, hence a section of the
logarithmic tangent bundle $T_X(-\log E)$. At a smooth point of
$E=\{x=0\}$ this bundle is generated by $x\partial_x$ and $\partial_y$;
at a crossing $E=\{xy=0\}$ it is generated by $x\partial_x$ and
$y\partial_y$. The logarithmic and ordinary saturations must be
distinguished: the ordinary foliation is obtained by dividing the lifted
field by its full divisorial zero factor as a section of $T_X$.
An exceptional component is invariant when the resulting foliation is
tangent to it at its general point, and dicritical when the foliation is
generically transverse to it. A dicritical component yields infinitely many
distinct separatrices on the original germ by taking the images of transverse
leaves. In this paper the resolution of the surface is not, unless stated
otherwise, a reduction of the singularities of the foliation in the sense of
Seidenberg \cite{Seidenberg68,Santos20}.

We use the smooth-branch convention for the GSV index as in
\cite[\S3 and Proposition~4]{Brunella97}, and the Camacho--Sad
normal-bundle index of \cite{CamachoSad82}. All indices below are
computed on the smooth resolution $X$. Suppose that
$C=\{x=0\}$ is a smooth invariant curve through $q=(0,0)$. An ordinary
saturated local generator and its dual one-form may then be written as
\[
 \xi=x a(x,y)\partial_x+b(x,y)\partial_y,\qquad
 \omega=b(x,y)\,dx-xa(x,y)\,dy.
\]
Here $b(0,y)$ is not identically zero. The G\'omez-Mont--Seade--Verjovsky
index of the foliation along this smooth branch and the Camacho--Sad index
are, with our conventions,
\begin{equation}\label{eq:prelim-indices}
 \GSV_q(\FF,C)=\ord_{y=0}b(0,y),\qquad
 \operatorname{CS}_q(\FF,C)=\operatorname{Res}_{y=0}
       \frac{a(0,y)}{b(0,y)}\,dy.
\end{equation}
Indeed, writing $\omega=g\,dx+x\eta$ with $g=b$ and
$\eta=-a\,dy$, the smooth-branch GSV index is the vanishing order of
$g|_C=b(0,y)$. The Camacho--Sad residue is obtained from
$-\eta/g=a\,dy/b$ along $C$. Thus the first index is a nonnegative
integer, whereas the second may be any complex number. If $C$ is compact,
restriction of $T_{\FF}$ to $C$ and the Camacho--Sad theorem give
\begin{equation}\label{eq:prelim-index-sums}
 \sum_{q\in C}\GSV_q(\FF,C)
   =\deg T_C-c_1(T_{\FF})\cdot C,\qquad
 \sum_{q\in C}\operatorname{CS}_q(\FF,C)=C^2.
\end{equation}
The first equality is the degree of the tangential morphism
$T_{\FF}|_C\to T_C$; the second is a residue formula for the normal
bundle $N_{C/X}$. For example, take
$\xi=\lambda x\partial_x+\mu y\partial_y$ at a crossing, with
$\lambda\mu\ne0$. On $x=0$ the tangential component is $\mu y$ and
the normal linearisation is $\lambda x$, so the GSV and Camacho--Sad
indices are respectively $1$ and $\lambda/\mu$. On $y=0$ they are
$1$ and $\mu/\lambda$.

For a smooth compact curve $C$ which is not invariant under $\FF$, we
use Brunella's tangency number $\operatorname{tang}(C,\FF)$
\cite[pp.~571--573, Lemma~2]{BrunellaENS97}. At $q\in C$, choose coordinates
$(x,y)$ with $C=\{x=0\}$ and an ordinary saturated generator
$\xi=A\partial_x+B\partial_y$. The local multiplicity in Brunella's
notation is $I(q,C,\FF)=\dim_{\CC}\OO_{X,q}/(x,A)$;
it is the vanishing order of $A(0,y)$, including when $q$ is a
singularity of $\FF$. The total number is the degree of the tangency
divisor, equivalently the zero divisor of the nonzero morphism
$T_{\FF}|_C\to N_{C/X}$, and thus
$\operatorname{tang}(C,\FF)=C^2-c_1(T_{\FF})\cdot C$.
It is distinct from the GSV index on an invariant curve.

Finally, at an isolated singularity $q$ of a foliation on a smooth surface,
the Baum--Bott residue associated with $c_1^2$ is the Grothendieck residue
\begin{equation}\label{eq:prelim-bb}
 \operatorname{BB}_{c_1^2}(\FF,q)
 =\operatorname{Res}_q\left\{
 \frac{(\operatorname{tr}D\xi)^2}{A B}\,dx\wedge dy\right\}.
\end{equation}
For a nondegenerate singularity with eigenvalues $\lambda,\mu\ne0$,
the linear terms $A=\lambda x$ and $B=\mu y$ give
$\operatorname{Res}_q(dx\wedge dy/(xy))=1$; hence the residue is
$(\lambda+\mu)^2/(\lambda\mu)=2+\lambda/\mu+\mu/\lambda$. If $X$ is compact, the Baum--Bott
theorem states that
\[
 \sum_{q\in\operatorname{Sing}(\FF)}
       \operatorname{BB}_{c_1^2}(\FF,q)=c_1(N_{\FF})^2,
 \qquad N_{\FF}=\det T_X\otimes T_{\FF}^{-1}.
\]
This formula is used only on smooth resolutions, where its local residues
are well defined \cite{BaumBott72,Brunella97}. In particular, it must not
be confused with the degree-one Camacho--Sad identity on an exceptional
curve.

\section{Lifting and exceptional valuations}\label{sec:val}

Let $(Y,p)$ be a normal complex surface germ with an isolated singularity and locally free tangent sheaf $T_Y=\operatorname{Der}_{\CC}(\OO_Y)$. Write $R=\OO_{Y,p}$ and choose a frame $(e_1,e_2)$. Let $\pi\colon X\to Y$ be a good resolution on which the frame lifts logarithmically (for instance the minimal good resolution), with compact reduced exceptional divisor $E=\sum_{i=1}^rE_i$. On this resolution the frame lifts to two holomorphic vector fields tangent to $E$; they are sections of $T_X(-\log E)$. This lifting statement is \cite[Fact~1, including the surface case]{BergnerGraf20}: the functorial resolution of a normal surface is its minimal good resolution, and its pushforward logarithmic tangent sheaf is reflexive. On a further arbitrary point blow-up, the rational transform of an individual vector field need not be holomorphic. All effectivity statements and Smith valuation bounds below refer to resolutions on which the tangent frame lifts holomorphically; in particular, they hold on the minimal good resolution. For other good resolutions we use only the persistence of generic transversality on strict transforms of its exceptional components.

For $v=fe_1+ge_2$ the ideal $J(v)=(f,g)$ is independent of the frame: a change of frame multiplies its coefficient column by a matrix in $\operatorname{GL}_2(R)$. In particular, a field is primitive if and only if $J(v)=R$, while $J(v)$ is $\mathfrak m$-primary exactly when the coefficient zero set is isolated at $p$. Saturation is understood in the sense of Section~\ref{sec:prelim}; it does not require a reduction of the singularities of the foliation.

At the generic point $\eta_i$ of $E_i$, let $A_i=\OO_{X,\eta_i}$, a discrete valuation ring with uniformiser $t$, fraction field $K_i$ and residue field $\CC(E_i)$. In logarithmic coordinates and frames, the lift of the tangent frame is represented by a full-rank matrix
\[
 M_i\colon A_i^2\longrightarrow A_i^2,
 \qquad (\widetilde e_1,\widetilde e_2)= (t\partial_t,\partial_z)M_i.
\]
The usual changes of the two frames multiply $M_i$ on the left and right by invertible matrices over $A_i$. Hence the two Smith exponents $0\leq\alpha_i\leq\beta_i$ are intrinsic to this morphism. Its determinant has valuation
\begin{equation}\label{eq:det}
 w_i=\ord_{E_i}(\det M_i)=\alpha_i+\beta_i.
\end{equation}
The frame also trivialises $\det T_Y$ in a neighbourhood of $p$. Since $Y$ is normal and $T_Y$ is free, $\omega_Y=(\det T_Y)^*$ is invertible; we may choose a nowhere-zero reflexive two-form $\nu$ on the smooth locus such that $\nu(e_1,e_2)=1$. At a general point of $E_i=\{t=0\}$ we can write
\[
 \pi^*\nu=t^{a_i}u(t,z)\,dt\wedge dz,\qquad u(0,z)\ne0,
\]
where $a_i$ is the discrepancy in $K_X=\pi^*K_Y+\sum_i a_iE_i$. Evaluating the two-form on $\widetilde e_1,\widetilde e_2$
gives $1=\pi^*\nu(\widetilde e_1,\widetilde e_2)$
on $X\setminus E$. Since
$\pi^*\nu(t\partial_t,\partial_z)=t^{a_i+1}u$,
this equality reads $1=t^{a_i+1}u\det M_i$.
Taking the $E_i$-valuation and using \eqref{eq:det} gives
\begin{equation}\label{eq:discrepancy}
 w_i=\alpha_i+\beta_i=-a_i-1\ge0,\qquad
 W=\sum_iw_iE_i=-\sum_i(a_i+1)E_i.
\end{equation}
The cycle $W$ is the determinant divisor of the logarithmic lifting map.

\begin{theorem}\label{thm:valuation}
For a nonzero field $v=f e_1+g e_2$, put $m_i=\min\{\ord_{E_i}(\pi^*f),\ord_{E_i}(\pi^*g)\}$ and let $z_i$ be the order of the logarithmic divisorial zero of $\widetilde v$ along $E_i$. There exist matrices $U_i,V_i\in\operatorname{GL}_2(A_i)$ such that
\[
 U_iM_iV_i=\begin{pmatrix}t^{\alpha_i}&0\\0&t^{\beta_i}\end{pmatrix}.
\]
Writing $(q_1,q_2)^{\mathsf t}=V_i^{-1}(\pi^*f,\pi^*g)^{\mathsf t}$, one has the exact identity
\begin{equation}\label{eq:exact}
 z_i=\min\{\alpha_i+\ord_{E_i}(q_1),\ \beta_i+\ord_{E_i}(q_2)\},
\end{equation}
and therefore
\begin{equation}\label{eq:bound}
 m_i+\alpha_i\leq z_i\leq m_i+\beta_i.
\end{equation}
If $\alpha_i<\beta_i$, equality on the left holds precisely when $\ord_{E_i}(q_1)=m_i$. Equivalently, the leading projective coefficient direction avoids the kernel of the reduction of $t^{-\alpha_i}M_i$ along $E_i$.
\end{theorem}
\begin{proof}
Smith normal form applies because $A_i$ is a DVR and $M_i$ has full rank over $K_i$. As $M_i=U_i^{-1}\operatorname{diag}(t^{\alpha_i},t^{\beta_i})V_i^{-1}$, multiplication by $U_i^{-1}$ does not change the minimum valuation of a vector's two coordinates; this yields \eqref{eq:exact}. The invertibility of $V_i$ gives $\min(\ord q_1,\ord q_2)=m_i$. The lower bound is immediate. At least one of $q_1,q_2$ has valuation $m_i$, so one entry in \eqref{eq:exact} has valuation at most $m_i+\beta_i$, giving the upper bound. If $\alpha_i<\beta_i$, a first coordinate of valuation $m_i$ gives equality on the left; if it has larger valuation, both entries have strictly larger valuation than $m_i+\alpha_i$. Reducing $t^{-\alpha_i}M_i$ modulo $t$ gives a rank-one matrix when the exponents differ. Its kernel is the vanishing locus of that first Smith coordinate, proving the last assertion.
\end{proof}

The coefficient valuation $m_i$ and the logarithmic zero order $z_i$ are
therefore different invariants. In particular, the leading coefficient
direction, not merely the divisorial order of the ideal $(f,g)$,
determines whether equality holds in \eqref{eq:bound}.

\section{Saturation and indices along the exceptional divisor}\label{sec:index}

For the remainder of the section assume that $J(v)$ is either the unit ideal or $\mm$-primary. Thus $v$ has no zero on $Y\setminus\{p\}$ and all divisorial zeros of its lift are exceptional. Let
\[
 Z_v=\sum_i z_iE_i
\]
be the zero divisor of $\widetilde v$ as a section of $T_X(-\log E)$. Denote by $D_v$ the \emph{reduced} sum of those exceptional components which are dicritical for the ordinary saturated foliation $\FF_v$. We reserve $E$ for the complete exceptional divisor and $D_v$ for its dicritical part.

\begin{lemma}\label{lem:ordinary}
The zero divisor of $\widetilde v$, considered as an ordinary holomorphic vector field, is exactly
\begin{equation}\label{eq:N}
 N_v=Z_v+D_v.
\end{equation}
Consequently, after dividing out its divisorial zeros, the ordinary saturated field defines a section of $T_X\otimes\OO_X(-N_v)$ with no divisorial zeros.
\end{lemma}
\begin{proof}
At a general point of $E_i=\{t=0\}$, the logarithmic lift has the form
\[
 \widetilde v=t^{z_i}\bigl(t A(t,z)\partial_t+B(t,z)\partial_z\bigr),
\]
with $A(0,z)$ and $B(0,z)$ not both identically zero. If $B(0,z)$ is nonzero, the ordinary order is $z_i$ and $E_i$ is invariant. Otherwise $B=tB_1$, $A(0,z)$ is nonzero, and the ordinary field has exactly one additional factor $t$; after division it is generically transverse to $E_i$. There are no other divisorial zeros, by the assumption on $J(v)$. This proves \eqref{eq:N}.
\end{proof}

For each $E_i$, let $g_i=g(E_i)$ and let $\delta_i^{\mathrm{inv}}$ (respectively, $\delta_i^{\mathrm{dic}}$) be the total number of intersection points, counted with their SNC multiplicities, with neighbouring invariant (respectively, dicritical) components. We use the GSV index on invariant components and the classical tangency number on dicritical components, as defined in Section~\ref{sec:prelim}.

\begin{theorem}\label{thm:indices}
Under the preceding assumptions, the following identities hold.

If $E_i$ is invariant, then
\begin{equation}\label{eq:gsv}
 \sum_{q\in E_i}\GSV_q(\FF_v,E_i)
 =2-2g_i-N_v\cdot E_i
 =(W-Z_v)\cdot E_i+\delta_i^{\mathrm{inv}}.
\end{equation}
If $E_i$ is dicritical, the tangency number of the ordinary saturated foliation with $E_i$ is
\begin{equation}\label{eq:normal}
 \operatorname{tang}(E_i,\FF_v)
 =E_i^2-N_v\cdot E_i
 =-Z_v\cdot E_i-\delta_i^{\mathrm{dic}}.
\end{equation}
All terms on the left are nonnegative integers. In particular, if $E_i$ is invariant then $(W-Z_v)\cdot E_i\geq0$.
\end{theorem}
\begin{proof}
Let $\xi$ denote a local generator of the ordinary saturated line. By Lemma~\ref{lem:ordinary}, $\widetilde v$ has divisor $N_v$ in that line, so the line bundle is identified with $\OO_X(N_v)$ over a neighbourhood of the exceptional fibre. Thus its restriction to an invariant component has a nonzero morphism
\[
 \OO_X(N_v)|_{E_i}\longrightarrow T_{E_i},
\]
which is the tangential component of $\xi$. Its zero divisor has degree
\[
 \deg T_{E_i}-N_v\cdot E_i=2-2g_i-N_v\cdot E_i.
\]
To identify its local multiplicities, choose coordinates $(x,y)$ in which $E_i=\{x=0\}$. Invariance means
\[
 \xi=xA(x,y)\partial_x+B(x,y)\partial_y,
 \qquad \omega=B(x,y)\,dx-xA(x,y)\,dy.
\]
The local equation of the branch is $f=x$, and, as in \eqref{eq:prelim-indices}, $\omega=g\,df+f\eta$ with $g=B$ and $\eta=-A\,dy$. At a point where the tangential coefficient is not identically zero along the branch, the smooth-branch GSV index is precisely
\[
 \GSV_q(\FF_v,E_i)=\ord_q\bigl(B(0,y)\bigr).
\]
All such orders are nonnegative. Since the restriction is a nonzero section on the compact smooth curve $E_i$, their sum is its degree, in accordance with the first identity of \eqref{eq:prelim-index-sums}, proving the first equality of \eqref{eq:gsv}. In particular this formula counts a singular point at a crossing if the tangential component vanishes there; it does not require a reduced singularity.

To rewrite the degree, put $E_{\ne i}=E-E_i$. Adjunction gives
\[
 2-2g_i=-(K_X+E_i)\cdot E_i,
 \qquad W\cdot E_i=-(K_X+E_i+E_{\ne i})\cdot E_i.
\]
Since $E_i$ is invariant and $D_v$ contains only dicritical components,
$D_v\cdot E_i=\delta_i^{\rm dic}$, whereas
$E_{\ne i}\cdot E_i=\delta_i^{\rm inv}+\delta_i^{\rm dic}$. Inserting $N_v=Z_v+D_v$ gives
\[
 \begin{split}
 2-2g_i-N_v\cdot E_i
 &= (W-Z_v)\cdot E_i
    +\delta_i^{\rm inv}+\delta_i^{\rm dic}-D_v\cdot E_i\\
 &= (W-Z_v)\cdot E_i+\delta_i^{\rm inv}.
 \end{split}
\]
For the residual tangential zeros, use the logarithmic rather than the
ordinary line. At a crossing $E_i=\{x=0\}$, $E_j=\{y=0\}$, the
logarithmic tangent bundle is generated by $x\partial_x$ and
$y\partial_y$. A logarithmic field whose restriction is tangent to
$E_i$ has tangential coefficient divisible by $y$. This holds whether
or not $E_j$ is invariant under the ordinary saturated foliation.
Removing the compulsory crossing factors therefore gives the morphism
\[
 \OO_X(Z_v)|_{E_i}\longrightarrow
 T_{E_i}(-E_{\ne i}|_{E_i}).
\]
Because the field is not identically zero tangentially at the generic point of $E_i$, this is a nonzero map of line bundles. Its zero divisor has degree
\[
 2-2g_i-E_{\ne i}\cdot E_i-Z_v\cdot E_i
 =-(K_X+E)\cdot E_i-Z_v\cdot E_i
 =(W-Z_v)\cdot E_i\ge0.
\]
The degree is nonnegative because the morphism does not vanish
generically on $E_i$. In the ordinary GSV count, a crossing with an
invariant neighbouring curve contributes a compulsory zero, while a
crossing with a dicritical neighbour does not; the two effects are
exactly accounted for by $D_v\cdot E_i$ in \eqref{eq:gsv}.

If $E_i$ is dicritical, restriction of the ordinary saturated line has nonzero normal projection
\[
 \OO_X(N_v)|_{E_i}\longrightarrow N_{E_i/X}=\OO_{E_i}(E_i).
\]
The corresponding section of $N_{E_i/X}\otimes\OO_X(-N_v)|_{E_i}$ has degree
$E_i^2-N_v\cdot E_i$. As $E_i$ occurs in $D_v$ with coefficient one, one has
$D_v\cdot E_i=E_i^2+\delta_i^{\rm dic}$, proving \eqref{eq:normal}. The normal projection defines the classical tangency divisor on $E_i$; its local multiplicity counts a zero at a singular point of the foliation as well as a tangency at a regular point.
\end{proof}

\begin{remark}\label{rem:cs-gsv}
For the local generator $\xi=bx\partial_x-ay\partial_y$ at a crossing of invariant curves $\{x=0\}$ and $\{y=0\}$, with $ab\ne0$, the GSV index on either branch is $1$, because the tangential coefficient has a simple zero. The Camacho--Sad indices, on the other hand, are $-b/a$ and $-a/b$, and their product is $1$. These are different local invariants. Their conflation would invalidate both the corner correction in \eqref{eq:gsv} and the use of Camacho--Sad summation in examples.
\end{remark}

\begin{remark}\label{rem:corners}
The correction $\delta_i^{\mathrm{inv}}$, rather than the total number of crossings, is essential. In local coordinates $E_x=\{x=0\}$ and $E_y=\{y=0\}$, consider the logarithmic fields
\[
 u_1=x\partial_x+y\partial_y,\qquad
 u_2=xy\partial_x+y\partial_y=y(x\partial_x+\partial_y).
\]
For $u_1$, both components are invariant. Its log-saturated tangential coefficient on $E_x$ is a unit, but the ordinary saturated GSV index at the crossing is $1$. For $u_2$, the ordinary saturation is $x\partial_x+\partial_y$; $E_x$ is invariant and $E_y$ is dicritical, and the GSV index on $E_x$ at the crossing is $0$. Thus replacing $\delta_i^{\mathrm{inv}}$ by the total number of crossings gives a false formula. These are local models for the index correction, not resolutions of singular surfaces.
\end{remark}

An elementary calculation also shows why a zero scheme cannot replace the foliated information. On the smooth plane, the fields $x\partial_x+y\partial_y$ and $x\partial_x-y\partial_y$ have the same coefficient ideal $(x,y)$ of colength one. In the blow-up chart $x=u$, $y=uw$, their saturated transforms are, respectively,
\[
 \partial_u \qquad\text{and}\qquad
 u\partial_u-2w\partial_w.
\]
The first is generically transverse to the exceptional curve $u=0$, whereas the second leaves it invariant. Indeed, the \emph{unsaturated} lift of the radial field is $u\partial_u$, and removal of its divisorial factor $u$ is necessary to obtain $\partial_u$. Both fields have the same coefficient valuations at the origin. This example explains why our local criteria concern primitive generators of a free tangent sheaf rather than numerical data from an arbitrary isolated zero.

\section{Dicriticality of tangent-field generators}\label{sec:primitive-dicritical}\label{sec:exceptional-correction}

Write $A=\sum_i\alpha_iE_i$ and $B=\sum_i\beta_iE_i$, so $W=A+B$.

\begin{lemma}\label{lem:nonzero}
If $p$ is singular, then $W\ne0$ and hence $B\ne0$.
\end{lemma}
\begin{proof}
Let $U$ be a sufficiently small neighbourhood of the compact SNC divisor $E$ in $X$. The logarithmic residue sequence is
\begin{equation}\label{eq:residue-seq}
 0\longrightarrow\Omega_X^1\longrightarrow\Omega_X^1(\log E)
 \xrightarrow{\operatorname{res}}\bigoplus_i\OO_{E_i}\longrightarrow0.
\end{equation}
Every $E_i$ is compact and connected, so a holomorphic section of $\OO_{E_i}$ is constant. If a logarithmic one-form $\eta$ on $U$ has residue vector $(r_i)$, exactness gives
\[
 \sum_i r_i c_1(\OO_X(E_i))=0\quad\hbox{in }H^1(U,\Omega_X^1).
\]
Indeed, in local equations $t_i=0$ the connecting morphism is represented by the Cech cocycle $d\log(t_{i,\alpha}/t_{i,\beta})$, the first Chern class of $E_i$. Restricting the equality to each compact component $E_j$ and taking degree yields
\begin{equation}\label{eq:res-matrix}
 \sum_i r_i(E_i\cdot E_j)=0\qquad(1\le j\le r).
\end{equation}
The exceptional intersection matrix is negative definite
\cite[Chapter~II]{Laufer71}; this is the intersection-theoretic
form of the contractibility criterion \cite{Grauert62}. In particular,
it is invertible. Thus $r_i=0$ for every $i$, and exactness
of \eqref{eq:residue-seq} then shows that every untwisted logarithmic
one-form on $U$ is holomorphic.

Suppose now that $W=0$. The determinant of the two lifted frame fields is a unit section of $\det T_X(-\log E)$ at every point of $U$: away from $E$ it is the pullback of a frame determinant, and its divisor on $U$ is precisely $W$. Hence $\widetilde e_1,\widetilde e_2$ form a frame of $T_X(-\log E)$ on $U$. Their dual basis consists of untwisted logarithmic one-forms, so both are holomorphic by the preceding paragraph. At a smooth point of $E$, however, the values of the ordinary vector fields $\widetilde e_1,\widetilde e_2$ lie in the one-dimensional subspace $T E$. Two holomorphic one-forms cannot pair with these two vectors to give the $2\times2$ identity matrix. This contradiction proves $W\ne0$. Finally $\beta_i\ge\alpha_i\ge0$ and $W=A+B$, so $B\ne0$.
\end{proof}

Suppose the hypotheses of Section~\ref{sec:index} hold and every exceptional component is invariant. Write \(M=(E_i\cdot E_j)\) for the negative-definite intersection
matrix. Let
\[
 b_i=\sum_{q\in E_i}\GSV_q(\FF_v,E_i)-\delta_i^{\mathrm{inv}}.
\]
By the second zero count in the proof of
Theorem~\ref{thm:indices}, $b_i$ is the degree of the residual
logarithmic tangential zero divisor. In particular, it is nonnegative;
the subtraction removes one compulsory GSV contribution at every
crossing with an invariant component.

\begin{proposition}\label{prop:correction-cycle}
There is a unique exceptional cycle \(B_v\) satisfying \(B_v\cdot E_i=-b_i\) for all \(i\). It is effective, integral, and equal to \(Z_v-W\). In particular, writing \(B_v=\sum_i d_i E_i\), one has
\[
 (d_i)_i=-M^{-1}(b_i)_i,\qquad a_i=d_i-z_i-1 .
\]
If the fibre is connected and some \(b_i\) is nonzero, then every \(d_i\) is positive. If at some component \(m_i<\alpha_i\), the foliation must instead have a dicritical component. The same holds if \(m_i<\beta_i\) at a component where the first Smith coordinate has order \(m_i\).
\end{proposition}

\begin{proof}
Because all components are invariant, $D_v=0$ and $N_v=Z_v$. For each $E_i$ the logarithmic tangent projection in the proof of Theorem~\ref{thm:indices} has a zero divisor of nonnegative degree, giving
\begin{equation}\label{eq:b-system}
 b_i=(W-Z_v)\cdot E_i\ge0.
\end{equation}
Set $B_v=Z_v-W=\sum_i d_iE_i$, so $B_v\cdot E_i=-b_i$. With $M=(E_i\cdot E_j)$ this is the linear system $M\mathbf d=-\mathbf b$. Negative definiteness makes $M$ invertible, so the solution is unique and is integral because it is already represented by $Z_v-W$.

Effectivity follows directly from negative definiteness. Write $B_v=B_v^+-B_v^-$ with disjoint effective supports. If $B_v^-\ne0$, then
\[
 B_v\cdot B_v^-=B_v^+\cdot B_v^--(B_v^-)^2>0,
\]
using nonnegative intersections between distinct exceptional curves and strict negativity of $(B_v^-)^2$. But every coefficient of $B_v^-$ is nonnegative and $B_v\cdot E_i=-b_i\le0$, giving $B_v\cdot B_v^-\le0$, a contradiction. Thus $B_v$ is effective. If a coefficient $d_i$ vanishes, the equality
$-b_i=B_v\cdot E_i=\sum_{j\ne i}d_j(E_j\cdot E_i)$ forces every neighbour of $E_i$ to have coefficient zero. The exceptional fibre is connected, so either $B_v=0$ or every coefficient is strictly positive.

The discrepancy formula \eqref{eq:discrepancy} says
$w_i=-a_i-1$. Therefore $d_i=z_i-w_i=z_i+a_i+1$, or
$a_i=d_i-z_i-1$. This also identifies the correction cycle with the
difference between the logarithmic zero divisor and the determinant
zero divisor, independently of the chosen tangent frame. If the foliation is non-dicritical, $z_i=w_i+d_i\ge\alpha_i+\beta_i$. The Smith bound gives $z_i\le m_i+\beta_i$. Hence $m_i\ge\alpha_i+d_i$. The two stated special cases follow immediately (the second also follows from the exact formula \eqref{eq:exact}).
\end{proof}

If all exceptional components are invariant, their residual tangential
zero counts satisfy $b_i\ge0$ by \eqref{eq:b-system}. The case
$b_i=0$ for every $i$ gives $B_v=0$; the proof of
Theorem~\ref{thm:primitive-dicritical} excludes it on a singular germ.
Thus a non-dicritical field on a singular surface must have a residual
tangential zero somewhere on the exceptional divisor. In particular,
when all exceptional crossings are nondegenerate, there must be a
singularity on the smooth part of the exceptional divisor.

Dicriticality persists on resolutions dominating the minimal good resolution.

\begin{theorem}\label{thm:primitive-dicritical}
Let \((Y,p)\) be a singular normal complex surface germ with locally free tangent sheaf. Every primitive local vector field has a dicritical component on the minimal good resolution. The strict transform of that component is dicritical on every good resolution.

More generally, let \(v\) have no zeros on the punctured germ, and suppose its saturated foliation leaves every exceptional component invariant on the minimal good resolution \(\pi:X\to Y\). Write \(Z_v\) for the divisor of zeros of its holomorphic logarithmic lift and \(d_i=z_i+a_i+1\), where \(z_i=\ord_{E_i}Z_v\) and \(a_i\) is the discrepancy of \(E_i\). Then every \(d_i\) is positive. At each exceptional component one has
\begin{equation}\label{eq:primitive-valuation}
 m_i\geq d_i+\alpha_i>0,
\end{equation}
where \(m_i\) is the order of the coefficient ideal and \(\alpha_i\) the smaller Smith exponent of the lifted frame. No analogous positivity claim is made here for the divisorial orders of meromorphic transforms at additional blow-up centres.
\end{theorem}

\begin{proof}
We first work on the minimal good resolution, where the frame lifts holomorphically to logarithmic fields. Let $e_1,e_2$ be a local tangent frame and let $\nu$ be the dual reflexive volume form, normalised by $\nu(e_1,e_2)=1$ on the smooth locus. Write $K_X=\pi^*K_Y+\sum_i a_iE_i$, and let $W=\sum_i(-a_i-1)E_i$ be the determinant divisor from \eqref{eq:discrepancy}. Every component is assumed invariant for the ordinary saturated foliation. Consequently, its logarithmic and ordinary saturated lines agree at the generic point of each $E_i$, and the divisor of the lifted field in that line is $Z_v=\sum_i z_iE_i$.

Restriction to $E_i$ produces a nonzero morphism
\[
 \OO_X(Z_v)|_{E_i}\longrightarrow
 T_{E_i}(-E_{\ne i}|_{E_i}),
\]
whose zero divisor has degree
\begin{equation}\label{eq:proof-bi}
 b_i=-(K_X+E)\cdot E_i-Z_v\cdot E_i
     =(W-Z_v)\cdot E_i\ge0.
\end{equation}
Put $B_v=Z_v-W=\sum_i d_iE_i$. By the negative-part argument in Proposition~\ref{prop:correction-cycle}, $B_v$ is effective, and either every $d_i$ is positive or all $d_i$ are zero. The latter case has to be excluded.

Suppose $B_v=0$. On the punctured germ set $\eta=\iota_v\nu$. It is a holomorphic one-form there, but after pullback it is initially only meromorphic along the exceptional fibre. In coordinates $E_i=\{t=0\}$ at a general point, write
\[
 \widetilde v=t^{z_i}(tP\partial_t+Q\partial_z),\qquad
 \pi^*\nu=t^{a_i}u\,dt\wedge dz,
\]
where $u$ is a unit and $Q(0,z)$ is not identically zero by invariance and divisorial saturation. A direct contraction gives
\begin{equation}\label{eq:contraction}
 \pi^*\eta
 =t^{z_i+a_i+1}u(P\,dz-Q\,dt/t)
 =t^{d_i}u(P\,dz-Q\,dt/t).
\end{equation}
At a double point with local exceptional equation $xy=0$, the saturated field is $xP\partial_x+yQ\partial_y$, the dual volume has the form $x^{a_i}y^{a_j}u\,dx\wedge dy$, and the same calculation gives
\[
 \pi^*\eta=x^{d_i}y^{d_j}u(P\,dy/y-Q\,dx/x).
\]
If all $d_i$ vanish, the displayed expressions show that
$\pi^*\eta$ has at most logarithmic poles at the generic points of
all exceptional components. The sheaf $\Omega_X^1(\log E)$ is locally
free, hence reflexive, so this membership extends across the finite
set of exceptional crossings and the remaining codimension-two points.
The residue along $E_i$ is $-u(0,z)Q(0,z)$; because $E_i$ is
invariant and the ordinary saturated field is nonzero at its generic
point, $Q(0,z)$ is not identically zero. This contradicts the
vanishing of exceptional residues established in Lemma~\ref{lem:nonzero}. Consequently $B_v\ne0$, and connectedness of the exceptional
fibre implies $d_i>0$ for every $i$.

The exact valuation bound of Theorem~\ref{thm:valuation} now gives, on the minimal good resolution,
\[
 z_i\le m_i+\beta_i,
 \qquad z_i=w_i+d_i=\alpha_i+\beta_i+d_i.
\]
Subtracting $\beta_i$ proves $m_i\ge\alpha_i+d_i>0$. If $v$ is primitive, then $m_i=0$ because one coefficient of $v$ is a unit at $p$. This contradiction shows that some exceptional component is dicritical.

Let $\sigma:X'\to X$ be a good resolution dominating the minimal good resolution. At the generic point of the strict transform of any component $E_j$ of $E$, the map $\sigma$ is an isomorphism. A foliation generically transverse to $E_j$ is therefore generically transverse to its strict transform. This proves the persistence assertion. We do not use the preceding effectivity or Smith-valuation argument at additional exceptional curves: the rational pullback of a vector field through a point that is not fixed by its flow may have a pole, so the holomorphic logarithmic lifting hypothesis would no longer be automatic.
\end{proof}

The dicritical component has an immediate consequence for the invariant curves through $p$.

\begin{corollary}\label{cor:separatrices}
Let $(Y,p)$ be a normal surface germ with locally free tangent sheaf. Then $Y$ is smooth at $p$ if and only if there exists a primitive derivation whose foliation has only finitely many invariant irreducible analytic curve germs through $p$. In fact, every primitive derivation on a singular germ has infinitely many such germs, whereas a primitive derivation on a smooth surface has exactly one.
\end{corollary}
\begin{proof}
Suppose first that $Y$ is smooth and $v$ is primitive. In local coordinates $v$ is nonzero at $p$, so the holomorphic flow-box theorem gives coordinates $(s,t)$ in which $v=\partial_s$. A germ of an irreducible analytic curve through $p$ tangent to this regular foliation must be the plaque $t=0$. There is therefore exactly one invariant irreducible curve germ.

Now let $p$ be singular. By Theorem~\ref{thm:primitive-dicritical}, the saturated transform of every primitive $v$ is generically transverse to some component $E_i$. Choose a Zariski-open subset $E_i^\circ$ of smooth points of $E_i$ disjoint from all other exceptional components and from the isolated singularities and tangencies of the lifted foliation. For each $q\in E_i^\circ$, the flow-box theorem gives a small analytic disc $\Delta_q\subset X$ tangent to the foliation, meeting $E_i$ transversely at $q$ and no other component of $E$.

The map $\pi|_{\Delta_q}$ is a nonconstant holomorphic map of germs from a smooth disc to $(Y,p)$. Embed $(Y,p)$ in a smooth ambient space and choose a coordinate whose pullback to $\Delta_q$ is nonzero and vanishes to finite order at $q$. On sufficiently small representatives this coordinate makes $\pi|_{\Delta_q}$ a finite, hence proper, map onto its image. Remmert's proper mapping theorem then shows that the image is an irreducible analytic curve germ through $p$. The image is invariant under $v$, because $\pi$ intertwines the foliations away from the exceptional fibre. Distinct points $q$ give distinct curve branches: a fixed irreducible curve germ has a unique strict transform and hence only one limiting intersection point with a given component at a transverse regular point. As $E_i^\circ$ contains infinitely many points, the surface germ admits infinitely many invariant curve germs. This proves both claims.
\end{proof}

\begin{theorem}\label{thm:generic}
Assume that $p$ is singular. For a fixed local tangent frame $(e_1,e_2)$ there is an exceptional component $E_j$ and a finite set $S\subset\CC$ such that, for every $\lambda\notin S$, the ordinary saturated foliation of $v_\lambda=e_1+\lambda e_2$ is dicritical along $E_j$. More precisely, for all $\lambda\notin S$,
\[
 Z_{v_\lambda}=A,\qquad W-Z_{v_\lambda}=B,
 \qquad B\cdot E_j<0.
\]
\end{theorem}
\begin{proof}
For each $E_i$, formula \eqref{eq:exact} with $(f,g)=(1,\lambda)$ gives $z_i(v_\lambda)=\alpha_i$ for all but at most one $\lambda$: if $\alpha_i<\beta_i$, the first Smith coordinate of $(1,\lambda)$ is a nonzero linear function of $\lambda$ over the residue field $\CC(E_i)$, which cannot vanish identically for two distinct values; if $\alpha_i=\beta_i$ the claim holds for every $\lambda$. There are finitely many exceptional components, so outside a finite set all equalities hold simultaneously. The divisor $B$ is nonzero and effective by Lemma~\ref{lem:nonzero}. Negative definiteness implies $B^2<0$, hence some $E_j$ has $B\cdot E_j<0$. If $E_j$ were invariant for the saturated foliation of such a $v_\lambda$, the last assertion of Theorem~\ref{thm:indices} would give $(W-Z_{v_\lambda})\cdot E_j=B\cdot E_j\geq0$, a contradiction. Thus $E_j$ is dicritical for every $\lambda$ outside the finite exceptional set.
\end{proof}

The exceptional correction has a natural interpretation in terms of the
classical fundamental cycle.  Let
\(F_\pi=\sum_i\varphi_i E_i\) denote Artin's fundamental cycle: it is
the smallest nonzero effective exceptional cycle with
\(F_\pi\cdot E_i\leq0\) for all \(i\)
\cite{Artin66,Laufer71,NagyNemethi19}.  Its coefficients are positive because the
exceptional divisor is connected.  This construction depends only on the
intersection matrix, rather than on a vector field.

\begin{corollary}\label{cor:fundamental-bound}
Assume that \((Y,p)\) is singular, that \(T_Y\) is locally free, and
that \(v=fe_1+ge_2\) has no zero on the punctured germ.  If the ordinary
saturated foliation of \(v\) leaves every component of the minimal good
exceptional divisor invariant, then
\begin{equation}\label{eq:fundamental-defect}
 B_v=Z_v-W\geq F_\pi,\qquad
 z_i\geq\alpha_i+\beta_i+\varphi_i,\qquad
 m_i\geq\alpha_i+\varphi_i
 \quad\text{for every }i.
\end{equation}
In particular, if \(m_i<\alpha_i+\varphi_i\) for any exceptional
component, the saturated foliation of \(v\) has a dicritical component.
\end{corollary}
\begin{proof}
Proposition~\ref{prop:correction-cycle} and the proof of
Theorem~\ref{thm:primitive-dicritical} show that \(B_v\) is a nonzero
effective exceptional cycle and that \(B_v\cdot E_i\leq0\) for every
\(i\).  The defining minimality of \(F_\pi\) gives
\(B_v\geq F_\pi\), coefficient by coefficient.  The relation
\(W=A+B=\sum_i(\alpha_i+\beta_i)E_i\) gives the second inequality.
Finally the Smith valuation estimate \eqref{eq:bound} yields
\(z_i\leq m_i+\beta_i\); comparison with the second inequality
gives \(m_i\geq\alpha_i+\varphi_i\).  The last assertion is its
contrapositive.
\end{proof}

The last inequality can be compared with the maximal-ideal cycle, which
is generally different from the fundamental cycle
\cite{TomariTomaru26}.  Let
\begin{equation}\label{eq:max-ideal-cycle}
 M_\pi=\sum_i\mu_i E_i,\qquad
 \mu_i=\min_{h\in\mm}\ord_{E_i}(\pi^*h),
\end{equation}
where \(\mm\) is the maximal ideal of \(\OO_{Y,p}\).  This is the
divisorial part of the pullback ideal \(\mm\OO_X\).  Since the exceptional
part of a principal divisor pulled back from \(\mm\) is an effective
anti-nef cycle, \(M_\pi\geq F_\pi\) coefficientwise.  The following
consequence applies to genuinely nonprimitive fields.

\section{General fields and the maximal-ideal filtration}\label{sec:generic-filtration}

The inequalities of Corollary~\ref{cor:fundamental-bound} apply to arbitrary coefficient pairs. For general pairs in a fixed power of the maximal ideal the lifting order can be calculated exactly. The calculation also identifies which leading coefficient directions are lost under the lifting matrix.

Fix generators $h_1,\ldots,h_N$ of $\mm$. For $k\ge1$ let $V_k\subset
\mm^k$ be the complex vector space spanned by all degree-$k$
monomials in the $h_a$, and put $V_0=\CC$. We call a pair $(f,g)\in V_k^2$
\emph{general} if it belongs to a nonempty Zariski-open subset of this
finite-dimensional parameter space. The condition may depend on the chosen
finite generating set, but the generic divisorial orders below do not.

\begin{theorem}\label{thm:generic-adic}
Let $(Y,p)$ be a singular normal complex surface germ with $T_Y$ free,
and let $\pi:X\to Y$ be its minimal good resolution. With the notation
$A=\sum_i\alpha_iE_i$, $B=\sum_i\beta_iE_i$ and
$M_\pi=\sum_i\mu_iE_i$ introduced above, for every integer $k\ge0$
a general field $v=f e_1+g e_2$ with $(f,g)\in V_k^2$ satisfies
\begin{equation}\label{eq:generic-adic-order}
 Z_v=A+kM_\pi,\qquad W-Z_v=B-kM_\pi.
\end{equation}
For $k\ge1$, its coefficient ideal $(f,g)$ is $\mm$-primary.
If $E_j$ satisfies
\begin{equation}\label{eq:fixed-component-criterion}
 (B-kM_\pi)\cdot E_j<0,
\end{equation}
then $E_j$ is dicritical for the saturated foliation of every field
in a possibly smaller nonempty Zariski-open subset of $V_k^2$.
The exceptional locus of non-dicritical fields is further constrained by the linear conditions in Corollary~\ref{cor:exceptional-linear-locus}.
\end{theorem}
\begin{proof}
Fix a component $E_i$ and work over the discrete valuation ring
$\OO_{X,\eta_i}$, as in Theorem~\ref{thm:valuation}. By definition of
$\mu_i$, one of the generators $h_a$ has order $\mu_i$; hence
$h_a^k\in V_k$ has order $k\mu_i$, while every element of $V_k$
has order at least $k\mu_i$. The same statement holds for $k=0$,
with the unit $1\in V_0$.

Choose a Smith decomposition of the logarithmic lifting matrix at
$E_i$. If $(q_1,q_2)^{\mathsf t}=V_i^{-1}(\pi^*f,\pi^*g)^{\mathsf t}$,
its first coordinate has the form $q_1=u\pi^*f+v\pi^*g$, where
$u,v$ belong to the DVR and at least one is a unit. The leading
coefficient
\[
 (f,g)\longmapsto
 \overline{t^{-k\mu_i}(u\pi^*f+v\pi^*g)}\in\CC(E_i)
\]
is a nonzero complex-linear map on $V_k^2$. Indeed, choosing the
coordinate with unit coefficient equal to $h_a^k$ and the other
coordinate equal to zero gives a nonzero image. The kernel is therefore
a proper complex-linear subspace. Off this kernel,
$\ord_{E_i}(q_1)=k\mu_i$; since the second Smith coordinate has
order at least $k\mu_i$, the exact valuation formula
\eqref{eq:exact} gives
\[
 z_i=\min\{\alpha_i+k\mu_i,\,\beta_i+\ord_{E_i}(q_2)\}
      =\alpha_i+k\mu_i.
\]
There are only finitely many exceptional components, so deleting their
finitely many proper linear exceptional loci gives the simultaneous
equality $Z_v=A+kM_\pi$. The determinant identity $W=A+B$
gives the second equality in \eqref{eq:generic-adic-order}.

For $k\ge1$, the general parameter pair is also a system of
parameters. One way to see this is to choose linear combinations
$\ell_1,\ell_2$ of the $h_a$ which form a system of parameters, by
generic Noether normalisation. The pair $(\ell_1^k,\ell_2^k)$ belongs
to $V_k^2$ and has an isolated common zero. For openness, form the projective tangent cone
$C_pY=\operatorname{Proj}\operatorname{gr}_\mm(R)$, a proper curve.
The initial degree-$k$ forms of a parameter pair such as
$(\ell_1^k,\ell_2^k)$ have no common projective zero. The incidence
condition that two degree-$k$ forms have a common zero on $C_pY$
is closed, by properness of $C_pY$ over $\CC$. Its complement is
therefore a nonempty Zariski-open subset of $V_k^2$, and any pair in
this complement generates an $\mm$-primary ideal. Intersecting this
open set with the finitely many Smith leading-term open conditions
proves the claim.

Suppose $E_j$ is invariant for the ordinary saturated foliation.
The logarithmic tangential zero count from
Theorem~\ref{thm:indices} gives
$(W-Z_v)\cdot E_j\ge0$. Equation~\eqref{eq:generic-adic-order}
contradicts this whenever \eqref{eq:fixed-component-criterion} holds.
Hence that component is dicritical on the claimed open set.

\end{proof}

The orders in Theorem~\ref{thm:generic-adic} do not determine which exceptional components are dicritical. The reduction of the leading logarithmic term distinguishes the two possibilities, and the generic distinction is independent of $k$.

\begin{theorem}\label{thm:generic-type}
In the setting of Theorem~\ref{thm:generic-adic}, there is a reduced exceptional divisor $D_{\mathrm{gen}}\subseteq E$, independent of $k\geq0$, with the following property: for every $k$, a general coefficient pair $(f,g)\in V_k^2$ defines an ordinary saturated foliation whose dicritical divisor is exactly $D_{\mathrm{gen}}$.

More explicitly, at the generic point of $E_i$, choose a local parameter $t$ and a logarithmic frame $(t\partial_t,\partial_z)$ and set
\[
 \overline{M}_i=\bigl(t^{-\alpha_i}M_i\bmod t\bigr)\in\operatorname{Mat}_2(\CC(E_i)),
\]
where the bar denotes reduction modulo $t$. Let $\ell_i\colon\CC(E_i)^2\to\CC(E_i)$ be projection to the coefficient of $\partial_z$. Then
\begin{equation}\label{eq:generic-type}
 E_i\subseteq D_{\mathrm{gen}}\quad\Longleftrightarrow\quad
 \ell_i\overline{M}_i=0.
\end{equation}
In particular every generically dicritical component has $\alpha_i<\beta_i$. If $p$ is singular, then $D_{\mathrm{gen}}\ne0$.

Write $s_i=-M_\pi\cdot E_i\geq0$ and let $\delta_{i,\mathrm{gen}}^{\mathrm{inv}}$ and $\delta_{i,\mathrm{gen}}^{\mathrm{dic}}$ count crossings with the invariant and dicritical components of $E-D_{\mathrm{gen}}$ and $D_{\mathrm{gen}}$, respectively. For general fields in $V_k^2$, the index formulae take the uniform form
\begin{equation}\label{eq:generic-index-variation}
 \begin{aligned}
 \sum_{q\in E_i}\GSV_q(\FF_{v_k},E_i)
 &= B\cdot E_i+ks_i+\delta_{i,\mathrm{gen}}^{\mathrm{inv}}
 && (E_i\not\subset D_{\mathrm{gen}}),\\
 \operatorname{tang}(E_i,\FF_{v_k})
 &= -A\cdot E_i+ks_i-\delta_{i,\mathrm{gen}}^{\mathrm{dic}}
 && (E_i\subset D_{\mathrm{gen}}).
 \end{aligned}
\end{equation}
Consequently $B\cdot E_i\geq0$ on a generically invariant component and $-A\cdot E_i\geq\delta_{i,\mathrm{gen}}^{\mathrm{dic}}$ on a generically dicritical one. If $s_i=0$, the corresponding index count is independent of $k$.
\end{theorem}
\begin{proof}
Fix $E_i$ and write $R_i=\OO_{X,\eta_i}$, a discrete valuation ring with residue field $K_i=\CC(E_i)$. The smaller Smith exponent $\alpha_i$ is the minimum valuation of the entries of $M_i$, so $\overline{M}_i$ is a nonzero matrix. Since $\alpha_i+\beta_i$ is the determinant order, $\overline{M}_i$ has rank one when $\alpha_i<\beta_i$, and rank two when $\alpha_i=\beta_i$.

Write $\mu_i=\ord_{E_i}(\mm\OO_X)$ and, for $(f,g)\in V_k^2$, consider its leading coefficient column
\[
 c_{i,k}(f,g)=\overline{t^{-k\mu_i}(\pi^*f,\pi^*g)^{\mathsf t}}\in K_i^2.
\]
This is a complex-linear map $V_k^2\to K_i^2$. Because one of the fixed generators $h_a$ of $\mm$ satisfies $\ord_{E_i}(h_a)=\mu_i$, the columns $(h_a^k,0)$ and $(0,h_a^k)$ have nonzero leading coefficients in the two respective coordinates. For $k=0$ use the constant columns $(1,0)$ and $(0,1)$. Thus the image of $c_{i,k}$ spans $K_i^2$ over $K_i$.

The coefficient column of the logarithmic lift is
\[
 M_i(\pi^*f,\pi^*g)^{\mathsf t}
 =t^{\alpha_i+k\mu_i}\bigl(\overline{M}_i c_{i,k}(f,g)+t(\cdots)\bigr).
\]
The condition $\overline{M}_i c_{i,k}(f,g)\ne0$ removes a proper complex-linear subspace of $V_k^2$. Outside it the logarithmic divisorial order is exactly $\alpha_i+k\mu_i$, as in Theorem~\ref{thm:generic-adic}. After division by this factor, the reduction of the logarithmic vector field is $\overline{M}_i c_{i,k}(f,g)$. The ordinary saturated foliation leaves $E_i$ invariant if its tangential component $\ell_i\overline{M}_i c_{i,k}(f,g)$ is nonzero; if this component vanishes and the logarithmic reduction is nonzero, the normal coefficient is a unit and the ordinary saturation removes one further factor of $t$, making the foliation transverse to $E_i$.

If $\ell_i\overline{M}_i\ne0$, the spanning property shows that $\ell_i\overline{M}_i c_{i,k}$ is a nonzero linear map, so its vanishing defines another proper linear subspace. Thus a general field is invariant along $E_i$. If $\ell_i\overline{M}_i=0$, all fields of logarithmic order $\alpha_i+k\mu_i$ are dicritical along $E_i$. This proves \eqref{eq:generic-type} for every $k$. Since the number of exceptional components is finite, all these requirements and, for $k\ge1$, the condition that $(f,g)$ be $\mm$-primary hold on a common nonempty Zariski-open subset of $V_k^2$.

The criterion is intrinsic. Although the matrix $\overline{M}_i$ depends on
the chosen local frames, the condition that its image be contained in
the logarithmic normal direction does not. Indeed, at the generic
point of $E_i$, the tangential projection $T_X(-\log E)|_{E_i}\to T_{E_i}$ is canonical; changing the parameter or logarithmic frame changes its coefficient by a unit. A change of tangent frame on $Y$ acts by an invertible matrix on the source of $M_i$, hence does not change the image of $\overline{M}_i$. If $\alpha_i=\beta_i$, $\overline{M}_i$ is invertible, and $\ell_i\overline{M}_i\ne0$; this proves the stated restriction on Smith exponents. The nonemptiness of $D_{\mathrm{gen}}$ when $p$ is singular follows by applying Theorem~\ref{thm:primitive-dicritical} to a general constant coefficient pair, which is primitive.

For fields in the common generic open set, Theorem~\ref{thm:generic-adic} gives $Z_{v_k}=A+kM_\pi$, while the ordinary dicritical divisor is the fixed $D_{\mathrm{gen}}$. The crossing numbers therefore do not depend on $k$. Substituting in \eqref{eq:gsv} and \eqref{eq:normal} gives \eqref{eq:generic-index-variation}. Nonnegativity of the residual tangential degree on invariant components gives $B\cdot E_i+ks_i\ge0$, and nonnegativity of the tangency number on dicritical components gives $-A\cdot E_i+ks_i-\delta_{i,\mathrm{gen}}^{\mathrm{dic}}\ge0$. Taking $k=0$ proves the numerical inequalities. The last assertion follows directly from $s_i=0$.
\end{proof}

The Smith exponents alone do not distinguish the two generic types. Over a discrete valuation ring, take the logarithmic frame in the order $(t\partial_t,\partial_z)$. The matrices $\operatorname{diag}(t,t^3)$ and $\left(\begin{smallmatrix}0&t^3\\t&0\end{smallmatrix}\right)$ have the same Smith exponents $(1,3)$. In the first case the reduction after division by $t$ has image in the logarithmic normal direction, so a generic coefficient column gives a dicritical foliation; in the second the reduction has image in the tangential direction, giving an invariant foliation. These are local models over a discrete valuation ring. The geometric datum in \eqref{eq:generic-type} is the leading image of the lifted frame, not merely the two determinant exponents.

\begin{corollary}\label{cor:exceptional-linear-locus}
For each $k\ge0$, let $\mathcal N_k\subseteq V_k^2$ consist of the nonzero coefficient pairs whose induced foliation has no dicritical exceptional component; for $k\ge1$ restrict to pairs with $\mm$-primary coefficient ideal. Then
\begin{equation}\label{eq:exceptional-linear-locus}
 \mathcal N_k\subseteq\bigcap_{E_i\subset D_{\mathrm{gen}}}
 \ker\bigl(\overline{M}_i\circ c_{i,k}\bigr)\subsetneq V_k^2.
\end{equation}
In particular, every such field has a strict jump of the logarithmic lifting order,
\[
 z_i(v)>\alpha_i+k\mu_i\qquad(E_i\subset D_{\mathrm{gen}}).
\]
The right-hand side of \eqref{eq:exceptional-linear-locus} is a proper complex-linear subspace determined by the leading coefficient maps along the generic dicritical components.
\end{corollary}
\begin{proof}
If $\overline{M}_i c_{i,k}(f,g)\ne0$ at a component of $D_{\mathrm{gen}}$, then the logarithmic order is $\alpha_i+k\mu_i$ and, by \eqref{eq:generic-type}, the tangential component of its logarithmic reduction is zero. Ordinary saturation therefore makes the component dicritical. A non-dicritical field must satisfy $\overline{M}_i c_{i,k}(f,g)=0$ for every $E_i\subset D_{\mathrm{gen}}$, and its divisorial order is consequently strictly larger than $\alpha_i+k\mu_i$. Each $\overline{M}_i c_{i,k}$ is a nonzero complex-linear map by the spanning argument in the preceding proof. Since $D_{\mathrm{gen}}\ne0$ at a singular point, their common kernel is proper.
\end{proof}

\begin{corollary}\label{cor:generic-maximal-ideal}
Let $d_j=-B\cdot E_j$ and $s_j=-M_\pi\cdot E_j\ge0$, with $M_\pi$ defined in \eqref{eq:max-ideal-cycle}.
If $d_j>0$ for some $j$, then $E_j$ is dicritical for a general field
in $V_k^2$ for \emph{every} integer $k\ge0$. In particular, the inequality
$d_j>ks_j$ remains a direct numerical test at a chosen order, but it is
not needed for persistence across the maximal-ideal filtration.
For an arbitrary field with no zeros off $p$ whose saturated foliation is non-dicritical, the individual valuation bound of Corollary~\ref{cor:fundamental-bound} remains applicable.
\end{corollary}
\begin{proof}
A general element of $\mm$ has exceptional divisor $M_\pi$:
a suitable linear combination of its fixed generators attains the
minimum valuation on all components simultaneously. Intersecting its
principal divisor with $E_j$ therefore gives
$M_\pi\cdot E_j\le0$, so $s_j\ge0$.
By the definitions, $(B-kM_\pi)\cdot E_j=-d_j+ks_j$.
If $d_j>0$, then $B\cdot E_j<0$, so Theorem~\ref{thm:indices}
excludes invariance for a general member of $V_0^2$.
Theorem~\ref{thm:generic-type} propagates dicriticity of this same
component to general fields in $V_k^2$ for all $k\ge0$.
The inequality $d_j>ks_j$ also implies dicriticity directly from
Theorem~\ref{thm:generic-adic}.
For non-dicritical fields outside the general set, the inequality in Corollary~\ref{cor:fundamental-bound} still applies.
\end{proof}

For special coefficient directions the logarithmic order can be strictly
larger than the generic order $\alpha_i+k\mu_i$. Corollary~\ref{cor:exceptional-linear-locus}
identifies the leading-term cancellations that must occur if a field in
$V_k^2$ is non-dicritical; the arbitrary-field inequality
\eqref{eq:fundamental-defect} then gives an independent numerical test.
The linear systems $V_k$ depend on the chosen generators of $\mm$,
but the generic orders and the divisor $D_{\mathrm{gen}}$ do not.

\section{Nonprimitive fields and exceptional indices}\label{sec:hamiltonian}
\subsection{First integrals}

Local first integrals are also compatible with a singular underlying surface. A first integral with normal crossings on the resolution gives a construction in which the indices are read directly from the exceptional multiplicities.

\begin{proposition}\label{prop:hamiltonian}
Let \((Y,p)\) be an isolated normal Gorenstein complex algebraic surface singularity, \(\pi:X\to Y\) a projective log resolution, and \(E=\sum_i E_i\) the reduced exceptional divisor. There are regular function germs \(f\in\mathfrak m_p^2\) and regular derivations \(v_f\) with the following properties. The derivation has no zeros on a punctured neighbourhood of \(p\), and its saturated foliation has first integral \(f\). On \(X\) every exceptional component is invariant, and all singularities over \(p\) are nondegenerate with negative rational eigenvalue quotient. If \(A_f=\sum_i n_iE_i\) is the exceptional part of \(\operatorname{div}(f\circ\pi)\), the indices on \(E_i\) are \(-n_j/n_i\) at a crossing with \(E_j\) and \(-1/n_i\) at any other singularity on \(E_i\). There are exactly \(-A_f\cdot E_i\) singularities of the second kind on \(E_i\).

Write $W_\pi=\pi^*K_Y-K_X-E$, an exceptional cycle which need not be effective on an arbitrary projective log resolution, and let $Z_{v_f,\pi}$ be the logarithmic zero divisor of the constructed field. Then
\begin{equation}\label{eq:ham-cycle-statement}
 Z_{v_f,\pi}-W_\pi=A_f.
\end{equation}
If $T_Y$ is locally free, the coefficient ideal of $v_f$ is $\mathfrak m_p$-primary. On resolutions where local tangent frames lift holomorphically and logarithmically, $W_\pi$ is the determinant divisor $W$ of Section~\ref{sec:exceptional-correction}.
\end{proposition}

\begin{proof}
Work on an affine representative and shrink it as necessary. Choose an integral exceptional cycle \(G>0\) with \(G\cdot E_i<0\) for every \(i\). For example, solve \(M\mathbf b=-\mathbf 1\) over \(\mathbb Q\), where \(M=(E_i\cdot E_j)\), and take \(G=\sum_i (N b_i)E_i\) with \(N\) a common denominator. The negative-part argument in Proposition~\ref{prop:correction-cycle} gives \(\mathbf b\geq 0\); since the right-hand side is strictly negative at every vertex and the exceptional fibre is connected, all its entries are positive. Put $\mathcal A=\OO_X(-G)$. Since $\deg(\mathcal A|_{E_i})=-G\cdot E_i>0$, the degree criterion for proper curves makes $\mathcal A$ ample on the possibly nonreduced exceptional fibre \cite[Tag~0B5Y]{StacksProject}. By \cite[Tag~0D2N]{StacksProject}, $\mathcal A$ is $\pi$-ample after shrinking the affine base. For $k\gg0$, relative Serre vanishing gives $R^1\pi_*(\mathcal A^{\otimes k}(-E))=0$ \cite[Tag~02O1]{StacksProject}. Since the base is affine, the restriction sequence yields a surjection $H^0(X,\mathcal A^{\otimes k})\to H^0(E,\mathcal A^{\otimes k}|_E)$. Increase $k$ so that $\mathcal A^{\otimes k}|_E$ is very ample and $n_i=2k\operatorname{coeff}_{E_i}(G)\geq a_i+1$, where $K_X=\pi^*K_Y+\sum_i a_iE_i$. By Bertini, choose a section on $E$ avoiding the nodes and with only simple zeroes on the smooth components; choose a second such section avoiding in addition the finite zero set of the first. Their zero divisors on $E$ are reduced and disjoint.

Lift both sections to $H^0(X,\mathcal A^{\otimes k})$. The inclusion $\OO_X(-kG)\subset\OO_X$ and $\pi_*\OO_X=\OO_Y$ give regular functions $f_1,f_2$ on the affine base. Write $g_i=\operatorname{coeff}_{E_i}(G)$. At a smooth point of $E_i=\{x=0\}$ we have $f_j\circ\pi=x^{kg_i}u_j$, where $u_j|_{E_i}$ is the chosen restricted section; at a node of $E_i$ and $E_\ell$ the expression is $x^{kg_i}y^{kg_\ell}u_j$ with $u_j$ a unit. Thus the exceptional divisor of each $f_j\circ\pi$ is exactly $kG$. A simple zero of $u_j|_{E_i}$ gives a smooth residual branch transverse to $E_i$, since its tangential derivative is nonzero. The residual divisors of $f_1$ and $f_2$ are disjoint near $E$, and properness allows the base to be shrunk so that no other branches intervene. Both $f_j$ vanish at $p$. Hence $f=f_1f_2\in\mathfrak m_p^2$ has $\operatorname{div}(f\circ\pi)=A_f+S$ near $E$, with $A_f=2kG$ and $S$ reduced, meeting $E$ simply at distinct smooth points. Principalness gives $S\cdot E_i=-A_f\cdot E_i$.

Choose a local generator \(\omega\) of the canonical sheaf of \(Y\). On the smooth locus define \(v_f\) by \(\iota_{v_f}\omega=df\). It extends uniquely as a regular derivation: the values on regular functions extend across \(p\) by normality, and the derivation rule persists. Furthermore \(df(v_f)=0\). Near a smooth point of the support of \(A_f+S\), the pulled-back function can be written after absorbing a unit as \(x^a\); near a crossing it can be written as \(x^a y^b\), with \(a,b>0\). At a smooth point of a component, a unit factor can be absorbed into
the transverse coordinate by a local holomorphic $a$-th root;
at a double point the units do not affect the saturated kernel.
Away from the exceptional fibre the strict transform $S$ is reduced,
so near $S\setminus E$ the function $f\circ\pi$ is a local coordinate.
At points outside $S\cup E$, critical points cannot accumulate onto
$E$: the displayed normal-crossings expressions show that the
differential is nonzero there. Properness permits shrinking the base
until $df$ has no zero on the punctured germ. The saturated kernels of the differentials are respectively generated by $\partial_y$ and $bx\partial_x-ay\partial_y$. At a crossing $E_i\cap E_j$, the unsaturated lifted field is a unit times
$x^{n_i-a_i-1}y^{n_j-a_j-1}(n_jx\partial_x-n_iy\partial_y)$; at a smooth point of $E_i$ the exceptional factor is $x^{n_i-a_i-1}$. The choice of $k$ makes these factors holomorphic, so the specially constructed derivation lifts logarithmically to this resolution even though a general tangent frame need not do so. In the second case $x=0$ has tangential eigenvalue $-a$ and normal eigenvalue $b$. The Camacho--Sad index with our convention is $-b/a$. Thus at an intersection of $E_i$ with $E_j$ it equals $-n_j/n_i$, and at a transverse point of the reduced strict transform $S$ it is $-1/n_i$. All these indices are negative rationals and the two eigenvalues are nonzero. Along an exceptional component away from these crossings the first integral is a power of a transverse coordinate times a unit; the saturated foliation is regular and tangent to that component. In particular, no saddle-node is introduced in the displayed resolution.

The number of intersections with $S$ follows directly from the principal divisor. The principal divisor $A_f+S$ has zero intersection with every compact exceptional component. Therefore
\begin{equation}\label{eq:section-count}
 \#(S\cap E_i)=S\cdot E_i=-A_f\cdot E_i,
\end{equation}
where all intersections are transverse and avoid the crossings of $E$. The GSV sum along $E_i$ is consequently
\[
 \sum_{q\in E_i}\GSV_q(\FF_{v_f},E_i)
 =(E-E_i)\cdot E_i+S\cdot E_i
 =(E-E_i)\cdot E_i-A_f\cdot E_i,
\]
which is the same tangential-degree computation as in the proof of Theorem~\ref{thm:indices}; on a resolution admitting holomorphic logarithmic lifts of the tangent frame it also agrees with \eqref{eq:gsv}.

The preceding normal-crossings calculation also proves that the
singularities of the saturated foliation over $p$ are precisely those
at the crossings of components of $A_f+S$. Its normal and tangential
eigenvalues are $b$ and $-a$, respectively, so neither eigenvalue
vanishes. The behaviour under further point blow-ups is recorded in
Proposition~\ref{prop:further-blowups}.

Suppose now that \(T_Y\) is locally free. In a local frame the canonical generator and the dual determinant frame differ by a unit on the smooth locus; the factor and its inverse extend by normality, so it is a unit on the germ. The coefficients of \(v_f\) are therefore, up to signs and a unit, the derivatives of \(f\) by the frame fields. They lie in \(\mathfrak m_p\) by the Leibniz rule and \(f\in\mathfrak m_p^2\). The zero at \(p\) is isolated by the preceding argument, hence the coefficient ideal is primary to \(\mathfrak m_p\). The local valuation can be computed explicitly. At a general point of $E_i=\{x=0\}$ write $f\circ\pi=x^{n_i}u(x,y)$, with $u$ a unit, and $\pi^*\omega=x^{a_i}h(x,y)\,dx\wedge dy$, with $h$ a unit. If $\widetilde v_f=P\partial_x+Q\partial_y$, then $\iota_{\widetilde v_f}(\pi^*\omega)=d(f\circ\pi)$ gives
\[
 \begin{split}
 P&=h^{-1}x^{n_i-a_i}\partial_yu,\\
 Q&=-h^{-1}x^{n_i-a_i-1}(n_i u+x\partial_xu).
 \end{split}
\]
Since $n_i u(0,y)\ne0$, the ordinary divisorial order is $n_i-a_i-1$. The saturated generator is tangent to $E_i$ and has nonzero tangential coefficient at its generic point; hence the logarithmic order is the same. The coefficient of $W_\pi=\pi^*K_Y-K_X-E$ at $E_i$ is $-a_i-1$, whether or not an entire tangent frame lifts holomorphically there. Consequently $z_i=n_i-a_i-1$, so that
\begin{equation}\label{eq:ham-cycle}
 Z_{v_f,\pi}=W_\pi+A_f.
\end{equation}
The correction $Z_{v_f,\pi}-W_\pi=A_f$ has intersection $-S\cdot E_i$ with $E_i$ by \eqref{eq:section-count}. Thus \eqref{eq:ham-cycle-statement} follows on the prescribed resolution without extending the determinant-divisor argument to an arbitrary blown-up tangent frame.
\end{proof}

\begin{corollary}\label{cor:integral-separatrices}
For the derivations constructed in Proposition~\ref{prop:hamiltonian}, the number of invariant irreducible curve germs through $p$ is finite. More precisely, every such germ is an irreducible component of the divisor $\{f=0\}$; on the chosen resolution the branches of the strict transform of that divisor meet the exceptional fibre at the $-A_f\cdot E_i$ points counted in \eqref{eq:section-count}.
\end{corollary}
\begin{proof}
Let $C$ be an invariant irreducible analytic curve germ through $p$, not contained in the singular point. On the dense open set of regular points of $C$ in the smooth locus of $Y$, the vector field $v_f$ is tangent to $C$. Since $v_f(f)=0$, the restriction of $f$ to this open set is constant. Continuity at $p$ forces this constant to be $f(p)=0$. Therefore $C$ is an irreducible component of the principal analytic divisor $f=0$. A nonzero function on a normal surface germ has only finitely many irreducible divisor components. Conversely, each reduced branch of $f=0$ is tangent to the kernel foliation of $df$ away from the critical point. This proves the assertion. The number of branches is given by \eqref{eq:section-count}, and the corresponding exceptional divisor is that of \eqref{eq:ham-cycle}.
\end{proof}

For the one-component analytic model of \cite[Proposition~5.4]{OAILipmanZariski}, the preceding construction, if the stated model exists, would give \(20k\) singularities of index \(-1/(2k)\) by taking \(G=C\): indeed, \(C^2=-10\) and \(A_f=2kC\). The coefficients of a nonprimitive field need not, however, come from a first integral. The following construction examines the freedom to prescribe exceptional indices using such coefficients. Unlike Proposition~\ref{prop:hamiltonian}, it uses the particular lifting and first-jet identities of the cited analytic model.

\subsection{Prescribed exceptional indices on an analytic model}\label{sec:prescribed-indices}

The following statement is conditional on the existence of the analytic
model described in \cite[Proposition~5.4]{OAILipmanZariski}, with its
stated lifting and first-jet identities. More precisely, we assume a
normal analytic surface germ $(Y,p)$ with a free tangent frame $(e_1,e_2)$
and a resolution having a single exceptional curve $C$ of genus $21$ and
self-intersection $-10$. Writing $L=\OO_C(-C)$, the assumptions include
$K_C\simeq L^4$ and the following first-jet data: in adapted
coordinates $C=\{t=0\}$, the ordinary lifting matrix is $t^2M$, with
\[
 M=\begin{pmatrix}r+tA_1+O(t^2)\\ B_1+tB_2+O(t^2)\end{pmatrix},
 \qquad \det M=t\delta,\quad\delta|_C\ne0,
\]
where $r=(r_1,r_2)=\xi(u_1,u_2)$, the $r_i$ form a basis of
$H^0(C,L)$, $\xi$ cuts out the six marked points $Q$, and the $u_i$
generate $L(-Q)$. The existence of such a model is not established
in this paper.

\begin{theorem}\label{thm:cs-realisation}
Under the preceding analytic-model hypotheses, write \(L=\OO_C(-C)\). For every sufficiently large integer \(m\), choose a section of \(L^{m-1}\) with simple zero divisor \(P\) disjoint from the six marked points of the construction. For every tuple
\[
  (\gamma_q)_{q\in P}\in\CC^P,\qquad \sum_{q\in P}\gamma_q=-10,
\]
there is a nonprimitive field \(v=f e_1+g e_2\), with \((f,g)\) primary to the maximal ideal of \(p\), whose saturated lift leaves \(C\) invariant and has precisely the points of \(P\) as singularities on \(C\). The tangential eigenvalues are nonzero and the Camacho--Sad index at \(q\) is \(\gamma_q\). The coefficient and logarithmic valuations are \(m_C=m\) and \(z_C=m+3\). In particular, the indices can all be chosen equal to \(-1/(m-1)\); for this choice the foliation is already reduced along \(C\) and has no saddle-nodes.
\end{theorem}

\begin{proof}
Let $I=\OO_X(-C)$ and let $Q$ be the divisor of the six marked
points. We use the normal-row and determinant identities in the
hypotheses, as taken from \cite[Proposition~2.3 and equations
(4.3), (4.11), (4.12), (5.5)]{OAILipmanZariski}. In particular the
frame is nonsingular away from $C$ after shrinking. The symbols for
the matrix rows refer to adapted local coordinates, whereas $r$ is
the global leading normal row.

For all \(n\) sufficiently large, holomorphic sections of \(I^n\) lift any given section of \(I^n/I^{n+1}=L^n\). Indeed, on every thickening \(jC\), the sheaf \(I^{n+1}/I^{n+1+j}\) has a filtration with quotients \(L^{n+1},\ldots,L^{n+j}\). These all have vanishing \(H^1\) by degree. The theorem on formal functions gives \(R^1\pi_*I^{n+1}=0\): the powers of \(I\) are cofinal near the compact fibre with the powers of the pullback of the maximal ideal, since the two ideals have the same radical there. The direct image vanishes away from \(p\), and its completion vanishes by the preceding computation. Take a sufficiently small Stein representative of \((Y,p)\); Leray and Cartan's theorem B give \(H^1(X,I^{n+1})=0\). The lifting assertion follows from the exact sequence for \(I^{n+1}\subset I^n\). The lifted functions descend to \(Y\) by normality.

At a generic point of \(C\), the rows \(r,B_1\) are dependent since \(\det M|_C=0\). Thus for \(h\in H^0(C,L^{m-1})\) with divisor \(P\), the column \(c_0=h(-r_2,r_1)^{\mathsf t}\) is killed by \(M|_C\), and this equality extends over \(Q\). Lift \(c_0\) to a pair of functions in \(I^m\), locally \(t^m(c_0+t c_1+\cdots)\). The resulting field vanishes to order at least \(m+3\) as an ordinary field. Its order-\((m+3)\) normal term lies in \(N_{C/X}\otimes L^{m+3}=L^{m+2}\); it is represented by \(r c_1+A_1c_0\), so it vanishes at \(Q\). The possible changes of this term from modifying the coefficient pair by \(I^{m+1}\) are the global sections in the image of
\[
 (L^{m+1})^{\oplus2}\xrightarrow{\,r\,} L^{m+2}(-Q).
\]
Since the sections \(r_1,r_2\) generate \(L(-Q)\), the kernel of this surjection is \(L^m(Q)\). It has vanishing \(H^1\) for \(m\) large. The lifting assertion for functions therefore permits a correction killing the normal term everywhere on \(C\).

Let \(s\) be the resulting tangential term, a section of \(T_C\otimes L^{m+3}\). For a constant local column \(e\), compare
\[
 \left.\det\bigl(M e,t^{-1}M(c_0+t c_1+\cdots)\bigr)\right|_{t=0}.
\]
The vanishing of the normal term gives \((re)s\) in local trivialisations. On the coefficient side this determinant is \(\delta|_C\det(e,c_0)=\delta|_C h(re)\). Taking an \(e\) with \(re\ne0\) off \(Q\) gives \(s=\delta|_C h\), and this identity extends across \(Q\). The multiplication is globally by a nowhere-zero identification \(L^{m-1}\simeq T_C\otimes L^{m+3}\). After division by \(I^{m+3}\), the ordinary field therefore leaves \(C\) invariant, with simple tangential zeros precisely at \(P\). In logarithmic frames the two Smith exponents at \(C\) are \(1,3\): the smallest entry has order \(1\) since the generic leading normal row is \(t r\), and the determinant has order \(4\). The coefficient order is \(m\), and the ordinary and logarithmic orders both equal \(m+3\) since the divided field is generically tangent and nonzero on \(C\).

We now vary the coefficient pair by sections in \(I^{m+2}\), keeping \(s\) fixed. Write the ordinary divided field locally as \(t a(t,z)\partial_t+b(t,z)\partial_z\). The change in \(a(0,z)\) is exactly \(r c_2\) for the order-\((m+2)\) pair \(c_2\). These changes are global sections of \(L^{m+3}(-Q)\), and all such sections occur, by the same kernel argument with \(m+2\) in place of \(m+1\). At \(q\in P\) the Camacho--Sad index is the residue of \(a(0,z)\,dz/b(0,z)\) for the smooth invariant branch \(t=0\). The differences of indices are consequently the residues at \(P\) of the differentials in \(H^0(C,K_C(P-Q))\): division by \(s\) identifies these differentials with \(H^0(C,L^{m+3}(-Q))\).

The only relation on the residue differences is that their sum vanish. To see this, first \(H^0(C,\OO_C(Q))=\CC\xi\). For any section \(\sigma\) in that space, both \(u_i\sigma\) belong to \(H^0(C,L)\) and therefore vanish at \(Q\), because \(r_1,r_2\) form a basis. Since the \(u_i\) have no common zero, \(\sigma\) also vanishes at \(Q\). Thus \(h^1(K_C(-Q))=1\) by Serre duality. The residue exact sequence for simple poles at \(P\), with the vanishing imposed at the disjoint divisor \(Q\), gives an image of codimension one in \(\CC^P\); indeed \(H^1(C,K_C(P-Q))=0\) as soon as \(\deg P>\deg Q\). The residue theorem identifies the image with the hyperplane of sum zero. The initial tuple of indices has sum \(C^2=-10\) by the Camacho--Sad index theorem. This proves the claimed realisation.

Finally, choose the first coefficient nonzero, as it already is at order \(m\). Perturb the second by an element of \(\mathfrak m^{m+3}\), which pulls back into \(I^{m+3}\) and hence leaves all the specified jets unchanged. There are finitely many height-one primes containing the first coefficient, and \(\mathfrak m^{m+3}\) is contained in none of them. In the complex vector space \(\mathfrak m^{m+3}\), the choices for which the perturbed second coefficient lies in one of those primes form a finite union of proper affine subspaces. Choose outside it. The coefficient ideal is then primary to \(\mathfrak m\), since both coefficients vanish at \(p\) and no height-one prime contains them both. On a smaller representative the field has no zeros off \(p\), so the division used above is the ordinary saturation and produces no other singularities near \(C\). With \(\deg P=10(m-1)\), the constant tuple in the statement sums to \(-10\), and every normal-to-tangent eigenvalue ratio is negative and nonzero. Hence every singularity on \(C\) is reduced and nondegenerate.
\end{proof}

\begin{corollary}\label{cor:no-go}
Assume that the singular surface germ with locally free tangent sheaf asserted in \cite[Proposition~5.4]{OAILipmanZariski} exists with the properties stated there. Then there is a nonprimitive derivation with an isolated coefficient zero whose saturated foliation has a non-dicritical reduction and only nondegenerate singularities with negative Camacho--Sad indices. Thus these dynamical properties alone would not imply smoothness, even if the tangent sheaf were locally free.
\end{corollary}
\begin{proof}
Apply Theorem~\ref{thm:cs-realisation} with the constant tuple $\gamma_q=-1/(m-1)$. The total sum is $-10$, since $\deg P=10(m-1)$. The reduced field has the nonzero tangential linear coefficient at each point of $P$ and a negative nonzero normal-to-tangential ratio. The surface is singular by the cited existence hypothesis. 
\end{proof}

For these fields, \(W=4C\), \(Z_v=(m+3)C\), and \(D_v=0\). Thus \((W-Z_v)\cdot C=10(m-1)\), as in Theorem~\ref{thm:indices}. The indices vary while the logarithmic valuation and the tangential zero divisor remain fixed.

The non-dicritical behaviour can also be checked after further blow-ups.
We use a cycle defined without lifting the tangent frame. On a good
resolution $\pi\colon X\to Y$ where every exceptional component is
invariant and all singularities along $E$ are nondegenerate, let
$b_i^{\mathrm{res}}$ be the number of singularities on the smooth part
$E_i^\circ=E_i\setminus\bigcup_{j\ne i}E_j$. Equivalently, it is the
degree of the zero divisor of the tangential morphism
$T_{\mathcal F}|_{E_i}\to T_{E_i}$ after subtracting the simple zero at
each invariant crossing. Negative definiteness determines a unique
rational exceptional cycle $B_{v,\pi}^{\mathrm{res}}$ by
\[
 B_{v,\pi}^{\mathrm{res}}\cdot E_i=-b_i^{\mathrm{res}}
 \qquad\text{for every }i.
\]
When a tangent frame lifts holomorphically and logarithmically to this
resolution, the restriction calculation of Theorem~\ref{thm:indices}
identifies $b_i^{\mathrm{res}}=(W-Z_v)\cdot E_i$; hence
$B_{v,\pi}^{\mathrm{res}}=Z_v-W$. We do not use this latter expression
on a resolution where the frame need not lift.

\begin{proposition}\label{prop:further-blowups}
Suppose that every exceptional component is invariant and all singularities
on $E$ are nondegenerate, with eigenvalue quotient in $\mathbb R_{<0}$.
Both properties persist after any sequence of point blow-ups over $p$.
For one blow-up $\sigma\colon X'\to X$ with exceptional curve $F$,
write $\pi'=\pi\circ\sigma$. The residual cycles satisfy
\[
 B_{v,\pi'}^{\mathrm{res}}=
 \begin{cases}
 \sigma^*B_{v,\pi}^{\mathrm{res}}+F,
  & \text{if the centre is singular on the smooth part of }E,\\
 \sigma^*B_{v,\pi}^{\mathrm{res}},
  & \text{otherwise}.
 \end{cases}
\]
\end{proposition}

\begin{proof}
A regular centre lies on a single branch, and flow-box coordinates give \(E=(y=0)\) and generator \(\partial_x\). In the chart \(x=u,\ y=uw\), its saturated transform is \(u\partial_u-w\partial_w\). The complementary chart \(x=uw,\ y=u\) gives the regular generator \(\partial_w\). Thus the new exceptional curve is invariant, with just one singularity, the crossing with the old branch. The new eigenvalue quotient is \(-1\).

At a singular centre use linear coordinates for the two eigenlines, with \(y=0\) the tangent line to one exceptional branch. Write the linear part as \(\lambda x\partial_x+\mu y\partial_y\). In the chart \(x=u,\ y=uw\) the transform has the form \(u(\lambda+O(u))\partial_u+((\mu-\lambda)w+O(u))\partial_w\). On the new exceptional line there is just one singularity in this chart, with eigenvalues \(\lambda,\mu-\lambda\); in the opposite chart the other singularity has eigenvalues \(\mu,\lambda-\mu\). Both quotients belong to \(\mathbb R_{<0}\), and the exceptional line is invariant. The strict transforms of the old invariant branches remain invariant and meet the new curve at the corresponding eigenlines. The transform has no exceptional common factor, since its component tangent to the new curve is not identically zero there. This computation also applies if one straightens the invariant branch before diagonalising the linear part.

The numbers $b_i^{\mathrm{res}}$ count precisely the nondegenerate
singularities away from the exceptional crossings. Blowing up a regular
point or a crossing leaves these numbers unchanged on the strict
transforms and gives $b_F^{\mathrm{res}}=0$. At a singular point on
$E_i^\circ$, one singularity is removed from that part of $E_i$ and
the new exceptional curve contains one singularity off the crossings:
thus $b_i^{\mathrm{res}}$ decreases by one and
$b_F^{\mathrm{res}}=1$. The projection formula gives
$\sigma^*B_{v,\pi}^{\mathrm{res}}\cdot E_i'= -b_i^{\mathrm{res}}$
for every strict transform and
$\sigma^*B_{v,\pi}^{\mathrm{res}}\cdot F=0$. In the last case,
$F\cdot E_i'=1$ and $F^2=-1$, so adding $F$ adjusts precisely the
two required intersection numbers. Uniqueness of the residual cycle
proves both transformation formulae. Iteration proves the assertion
for successive blow-ups.
\end{proof}

In the genus-\(21\) example the exhibited resolution is minimal: its only exceptional curve is not rational and cannot be an exceptional curve of the first point blow-down to a smooth surface. The minimal-resolution and factorisation theorems for normal analytic surface germs \cite[Chapter V]{Laufer71} therefore reduce every resolution which is an isomorphism off \(p\) to the calculation in Proposition~\ref{prop:further-blowups}. For the negative tuple in Theorem~\ref{thm:cs-realisation}, no further choice of resolution produces a dicritical component or a saddle-node. The residual cycle stays nonzero as well: the total number of singularities on the smooth parts of the exceptional components remains $10(m-1)$, although their distribution can change.

The prescribed-index construction is analytic. The corresponding algebraic statement requires a compatible analytic identification of the algebraisation in \cite[Proposition~6.2]{OAILipmanZariski}. Under that hypothesis the local coefficient scheme and the exceptional indices can be preserved.

\begin{corollary}\label{cor:algebraic-realisation}
Assume, in addition, that the affine algebraic surface asserted in \cite[Proposition 6.2]{OAILipmanZariski} exists, that its tangent sheaf is free, and that an analytic isomorphism identifies its germ with the analytic model in Theorem~\ref{thm:cs-realisation} together with the stated tangent-frame and resolution data. On that surface, one can realise each tuple in Theorem~\ref{thm:cs-realisation} by a regular algebraic vector field, after a principal affine shrinking about the singular point. On the minimal analytic resolution the field has the same coefficient valuation, the same divisorial saturation order, and the same first jet after saturation along \(C\) as the chosen holomorphic field. Its coefficient ideal has the same analytic germ. For the all-negative tuple, every algebraic resolution that is an isomorphism off the singular point is non-dicritical and has only nondegenerate singularities above that point.
\end{corollary}

\begin{proof}
Assume that the algebraic surface and tangent frame described in
\cite[Proposition~6.2]{OAILipmanZariski} exist as stated. The formal
identification of completed local rings invoked there is not used here
as a substitute for an analytic isomorphism. By the additional
hypothesis in this corollary, fix an analytic isomorphism compatible
with the tangent-frame and resolution data, and transport the analytic
resolution and field to $Y_{\mathrm{alg}}^{\mathrm{an}}$.

In the algebraic tangent frame let \(c=(c_1,c_2)^{\mathsf t}\) be the transported analytic coefficients and let \(J=(c_1,c_2)\). The change from the original frame is an invertible analytic matrix: algebraic tangent commutes with analytification, by the finite Jacobian presentation and flat base change. Hence \(J\) is unchanged. Both algebraic frame vectors lift to sections of \(I^2T_X\), since the old frame does and the entries of the change of basis pull back holomorphically; here \(I=\OO_X(-C)\).

Write \(\mathfrak m\) for the analytic maximal ideal at \(p\). Take \(s\) with \(\mathfrak m^s\subset J\), and choose
\[
 N\geq\max\{s+1,m+3\}.
\]
There are global regular functions \(a_1,a_2\) on \(Y_{\mathrm{alg}}\) with \(a_j-c_j\in\mathfrak m^N\): take degree-\((N-1)\) Taylor polynomials of convergent representatives in an affine embedding centred at \(p\). Since \(\mathfrak m^N\subset\mathfrak m J\), the new column is \((1+U)c\) for a matrix \(U\) with entries in \(\mathfrak m\). The matrix is invertible in the analytic local ring. Thus \((a_1,a_2)\) cuts out precisely the same analytic coefficient subscheme. The algebraic coefficient ideal is primary at \(p\) by faithful flatness of completion and the equality of the algebraic and analytic completed rings.

The pullback of \(\mathfrak m\) is contained in \(I\), so the two lifted fields differ by a section of \(I^{N+2}T_X\). The analytic lift has ordinary saturation order \(m+3\); the algebraic perturbation has the same order, since the generic leading nonzero term does not change. After saturation the difference lies in \(I^{N-m-1}T_X\subset I^2T_X\). In particular the two saturated fields agree on the first infinitesimal neighbourhood of \(C\). The invariance of \(C\), its tangential zero scheme, and the linearisation at every point of that scheme are preserved. Since the tangential eigenvalues at those points are nonzero, their Camacho--Sad indices are preserved as well. The old coefficient order \(m\) is unchanged, because the coefficient ideal itself is unchanged analytically.

The new field is regular algebraic on the affine model. Its coefficient zero locus has \(p\) as an isolated component, so a principal open neighbourhood of \(p\) excludes all other components. The frame there has no zero on the smooth locus, and neither does the new field away from \(p\). Finally, any algebraic resolution as in the statement analytifies near the fibre to a resolution dominating the minimal analytic one. Proposition~\ref{prop:further-blowups} applies to the all-negative tuple and gives the last assertion.
\end{proof}

\section{Baum--Bott residues on compact surfaces}\label{sec:compact-residue}

The exceptional correction also enters a global relation on a compact normal surface. In applying the Baum--Bott formula, we keep the residues of singularities on the smooth locus separate from the contributions of the exceptional divisor; no vanishing of these residues is assumed.

Let $Y$ be a compact normal complex surface with locally free tangent sheaf,
$\mathcal L$ a line bundle, and $v\in H^0(Y,T_Y\otimes\mathcal L)$ a
nonzero section with no divisorial zeroes. Its zeroes in $Y_{\mathrm{reg}}$
are then isolated and finite. Let $\pi:X\to Y$ be a good resolution which
is an isomorphism over $Y_{\mathrm{reg}}$, with reduced exceptional divisor
$E=\sum_i E_i$. Assume that on a neighbourhood of each exceptional fibre
\emph{every local frame of $T_Y$ lifts holomorphically to sections
of $T_X(-\log E)$};
thus the logarithmic lifting matrices and their determinant divisors
are defined as in Section~\ref{sec:val}. Suppose also that every $E_i$
is invariant under the ordinary saturated foliation $\mathcal F$ and
that all singularities of $\mathcal F$ on $E$ are nondegenerate. In particular, their two eigenvalues are nonzero. Define
\[
 \mathcal A=\mathcal L\otimes\omega_Y^{-1},\qquad
 H=\pi^*c_1(\mathcal A),\qquad
 B_v=Z_v-W=\sum_i d_iE_i.
\]
Here $Z_v$ is the logarithmic divisorial zero cycle, and
$W=\pi^*K_Y-K_X-E$ is the determinant divisor of the lifted tangent
frame. The displayed canonical-divisor expression follows from the
logarithmic determinant calculation \eqref{eq:discrepancy}
on the resolutions just specified; it patches over the exceptional
fibres. No lifting assertion is made for arbitrary further blow-ups. Put
\[
 \mathcal R=\sum_{x\in Y_{\mathrm{reg}}\cap\operatorname{Sing}(\mathcal F)}
               \operatorname{BB}_{c_1^2}(\mathcal F,x).
\]
At a degenerate isolated zero, the term is the usual Grothendieck
Baum--Bott residue, not an eigenvalue quotient. Let $P_i$ be the set of
singular points of $\mathcal F$ on the smooth part of $E_i$, and set
$P=\coprod_iP_i$. For $q\in P_i$ denote by $\gamma_q$ the
Camacho--Sad index on $E_i$. All intersections and squares below are
taken on the smooth compact surface $X$.

\begin{theorem}\label{thm:compact-residue}
Under these assumptions the following identity holds:
\begin{equation}\label{eq:compact-residue-general}
 H^2+B_v^2=\mathcal R+\sum_{q\in P}\gamma_q^{-1}.
\end{equation}
Suppose, in addition, that all Camacho--Sad indices of the exceptional
branches are negative real numbers. Orient each exceptional crossing
$q\in E_i\cap E_j$ once and write
$\operatorname{CS}(\mathcal F,E_i,q)=-p_q$, $p_q>0$; the index on
$E_j$ is then $-p_q^{-1}$. Write $\gamma_q=-r_q$, $r_q>0$, for
$q\in P_i$. Then
\begin{equation}\label{eq:compact-squares-r5}
 \mathcal R-H^2=
 \sum_{q:i\to j}
 \left(\sqrt{p_q}\,d_i-\frac{d_j}{\sqrt{p_q}}\right)^2
 +\sum_i\sum_{q\in P_i}
 \left(\sqrt{r_q}\,d_i-\frac{1}{\sqrt{r_q}}\right)^2.
\end{equation}
Consequently $\mathcal R\geq H^2$. If equality holds, and each
connected exceptional fibre contains at least one point of $P$, then
$d_i>0$ for every exceptional component, and the indices are necessarily
\begin{equation}\label{eq:rigid-indices-r5}
 \operatorname{CS}(\mathcal F,E_i,q)=-\frac{d_j}{d_i}
 \quad(q\in E_i\cap E_j),\qquad
 \gamma_q=-\frac1{d_i}\quad(q\in P_i).
\end{equation}
In particular, all exceptional indices are negative rational numbers
in this equality case.
\end{theorem}

\begin{proof}
Since $v$ has no divisorial zeroes off the exceptional set, its
saturated tangent line is
$T_{\mathcal F}=\pi^*\mathcal L^{-1}\otimes\OO_X(Z_v)$.
For the determinants of the logarithmic and ordinary normal sheaves,
respectively, this gives
\begin{equation}\label{eq:two-normals-r5}
 N_{\log}=c_1\bigl(\det T_X(-\log E)\otimes T_{\mathcal F}^{-1}\bigr)
 =H-B_v,\qquad N_{\mathrm{ord}}=N_{\log}+E.
\end{equation}
The first identity uses
$K_X+E=\pi^*K_Y-W$. Since $H\cdot E_i=0$,
we have $N_{\log}^2=H^2+B_v^2$.

We next calculate $N_{\log}\cdot E_i$ directly. The logarithmically
saturated field induces a nonzero tangential section on each invariant
$E_i$, whose zero divisor has degree $N_{\log}\cdot E_i$ by the
restriction calculation in Theorem~\ref{thm:indices}. At a nondegenerate
crossing of two invariant components this tangential section has nonzero
value in the logarithmic frame; at a nondegenerate singularity on the
smooth part of $E_i$ it has a simple zero. Hence
\begin{equation}\label{eq:boundary-degrees-r5}
 N_{\log}\cdot E_i=-B_v\cdot E_i=|P_i|,
 \qquad N_{\log}\cdot E=|P|.
\end{equation}
This also checks the correction at crossings independently of the
ordinary smooth-branch GSV count: the latter acquires one additional
zero at each invariant crossing, whereas the logarithmic tangential
section does not.

Let $c$ denote the number of crossings, with distinct intersection
points counted separately. At a nondegenerate singularity with eigenvalues $\lambda,\mu$,
formula~\eqref{eq:prelim-bb} gives
$(\lambda+\mu)^2/(\lambda\mu)=2+\gamma+\gamma^{-1}$,
where $\gamma=\lambda/\mu$. At a smooth-boundary singularity,
$\gamma$ is its Camacho--Sad index on the invariant exceptional
component; at a crossing the two branch indices are reciprocal.
Summing the
Camacho--Sad formula over all invariant components gives
\begin{equation}\label{eq:cs-summed-r5}
 \sum_i E_i^2=
 \sum_{q:i\to j}(\gamma_{i,q}+\gamma_{i,q}^{-1})
 +\sum_{q\in P}\gamma_q.
\end{equation}
The Baum--Bott theorem on the smooth compact surface $X$
\cite{BaumBott72} therefore yields
\[
 N_{\mathrm{ord}}^2=\mathcal R+\sum_iE_i^2
             +2(c+|P|)+\sum_{q\in P}\gamma_q^{-1}.
\]
On the other hand,
$E^2=\sum_iE_i^2+2c$, and
\eqref{eq:two-normals-r5}--\eqref{eq:boundary-degrees-r5}
show that
$N_{\mathrm{ord}}^2-N_{\log}^2=E^2+2|P|$.
Subtracting proves \eqref{eq:compact-residue-general}.
No term at a zero on $Y_{\mathrm{reg}}$ has been discarded.

For the second assertion, write $Q=-(E_i\cdot E_j)$ for the
negative of the exceptional intersection matrix and
$\mathbf d=(d_i)$. The negative definiteness of the exceptional
intersection matrix makes $Q$ positive definite, with nonpositive
off-diagonal entries. Equation~\eqref{eq:boundary-degrees-r5} says
\begin{equation}\label{eq:Qd-r5}
 Q\mathbf d=(|P_i|)_i,\qquad
 \mathbf d^{\mathsf t}Q\mathbf d=\sum_i d_i|P_i|=-B_v^2.
\end{equation}
By Camacho--Sad, the diagonal entries of $Q$ are the sums of the
positive numbers $p_q$, $p_q^{-1}$ and $r_q$ attached to the
corresponding branches; its off-diagonal entries are minus the numbers
of crossings. To see the quadratic identity, expand $\mathbf d^{\mathsf t}Q\mathbf d$
using the diagonal Camacho--Sad sums from
\eqref{eq:cs-summed-r5}: every crossing contributes
$p_qd_i^2+p_q^{-1}d_j^2-2d_id_j$ and every point of $P_i$
contributes $r_qd_i^2$. Thus
\begin{equation}\label{eq:energy-r5}
 \mathbf d^{\mathsf t}Q\mathbf d=
 \sum_{q:i\to j}
 \left(\sqrt{p_q}\,d_i-\frac{d_j}{\sqrt{p_q}}\right)^2
 +\sum_i\sum_{q\in P_i}r_qd_i^2.
\end{equation}
Expanding the final squares of \eqref{eq:compact-squares-r5}, and using \eqref{eq:energy-r5} and \eqref{eq:Qd-r5}, gives
\[
 \mathbf d^{\mathsf t}Q\mathbf d
 -2\sum_i d_i|P_i|+\sum_{q\in P}r_q^{-1}
 =-\mathbf d^{\mathsf t}Q\mathbf d+\sum_{q\in P}r_q^{-1}.
\]
The latter equals $\mathcal R-H^2$ by
\eqref{eq:compact-residue-general}, because
$B_v^2=-\mathbf d^{\mathsf t}Q\mathbf d$ and
$\gamma_q^{-1}=-r_q^{-1}$. This establishes the sum of squares.

For completeness, \eqref{eq:Qd-r5} implies $d_i\geq0$ without assuming
this in advance. Indeed, if the negative part $\mathbf d^-$ were
nonzero, multiplication of $Q\mathbf d\geq0$ by $\mathbf d^-$
would give a nonnegative number, whereas
$(\mathbf d^-)^{\mathsf t}Q\mathbf d<0$: the cross term between
the disjoint positive and negative supports is nonpositive, and the
quadratic form of $Q$ is strictly positive on $\mathbf d^-$. If
some $|P_i|$ is positive in a connected exceptional fibre, its
coefficients are all strictly positive: a zero coefficient would force
all adjacent coefficients to vanish by the nonpositive off-diagonal
entries, and then connectedness would force every $|P_i|$ to vanish.
When $\mathcal R=H^2$, each square in
\eqref{eq:compact-squares-r5} vanishes. The strict positivity just
proved allows us to divide by the $d_i$, giving
\eqref{eq:rigid-indices-r5}.
\end{proof}

\begin{corollary}\label{cor:nonresonant-smoothness-r5}
Let $Y$ be a normal projective surface with locally free tangent sheaf,
and suppose that a twisted field $v$ as above has no zero in
$Y_{\mathrm{reg}}$. Assume that $\mathcal A=\mathcal L\otimes\omega_Y^{-1}$
is nef and that on a good resolution the exceptional divisor is
invariant, with only nondegenerate singularities having negative real
Camacho--Sad indices. Suppose moreover that each nonempty connected
exceptional fibre contains a singularity away from the exceptional
crossings and that, whenever the exceptional divisor is nonempty,
at least one exceptional Camacho--Sad index is irrational. Then $Y$ is
smooth.
\end{corollary}
\begin{proof}
Nefness gives $H^2\geq0$, and the absence of zeroes on the regular
locus gives $\mathcal R=0$. If $E\ne\varnothing$, the sum of squares
in Theorem~\ref{thm:compact-residue} gives
$0\leq-H^2\leq0$. Thus $H^2=0$ and equality holds in every square.
Each exceptional fibre has a point of $P$ by hypothesis, so
\eqref{eq:rigid-indices-r5} makes all its Camacho--Sad indices
rational, contrary to the nonresonance assumption. Therefore $E$ is
empty, and normality implies that $Y$ is smooth.
\end{proof}

The rational ratios in \eqref{eq:rigid-indices-r5} are also
realised locally by the first-integral fields of
Proposition~\ref{prop:hamiltonian}. The compact smoothness criterion
uses an additional irrationality hypothesis. Under its negative-index
assumptions the exceptional foliation is reduced and has no saddle-nodes,
so the resulting smoothness conclusion is also related to the
more general generalised-curve criterion of
\cite[Corollary~1.9]{CorreaLourencoMachado26}. The contribution of
Theorem~\ref{thm:compact-residue} is the explicit correction-cycle
identity \eqref{eq:compact-residue-general} and its numerical
equality criterion \eqref{eq:rigid-indices-r5}. The underlying
comparison of logarithmic and ordinary normal sheaves belongs to the
residue theory developed in
\cite[Corollaries~1.6 and 1.7]{CorreaLourencoMachado26} and
\cite[Theorem~1.1]{CavalcanteCorreaLourencoShahsavaripour26}.

\enlargethispage{2\baselineskip}
\medskip
\noindent\textbf{Acknowledgments.}  
The author thanks Alana Cavalcante, Marcos Jardim, Fernando Louren\c{c}o,
Jos\'e Seade, Elaheh Shahsavaripour and Calum Spicer for useful
discussions on indices and vector fields on singular varieties.
He acknowledges partial support from the University of Bari Aldo Moro.
This work was completed during a visit to the Universidade Estadual de
Campinas (UNICAMP) as a visiting professor, supported by the
S\~ao Paulo Research Foundation (FAPESP), Visiting Researcher Programme,
grant 2026/19401-5.

\medskip
\noindent\textbf{AI usage statement.}
  The mathematical ideas and original results presented in this paper are due to the author. During the preparation of this paper, ChatGPT-6 Astra (OpenAI) was consulted for exploratory calculations and consistency checks, bibliographic enquiries, and the examination of possible counterexamples. It was also used to challenge proposed proof strategies and to scrutinise successive drafts from an adversarial referee’s perspective, with particular attention to gaps in the arguments and implicit assumptions. The author wrote the manuscript and takes full responsibility for the correctness of its mathematical content, including all statements, proofs and conclusions.

\end{document}